%% file: MMEOptimization.tex
\documentclass[preprint,3p,sort&compress]{elsarticle}
\journal{{\tt arXiv.org}}

\usepackage[dvips]{epsfig}
\usepackage{graphicx} 
\usepackage{pdfpages}
\usepackage{latexsym}
\usepackage{verbatim}
\usepackage{amsmath}
\usepackage{amsthm}
\usepackage{amssymb}
\usepackage[margin=10pt,font=small,labelfont=bf,
labelsep=endash]{caption}
\usepackage{subfig}
\usepackage{bbm}
\usepackage{pstool}
\usepackage[boxed]{algorithm2e}
\usepackage{mathtools}
\usepackage{fonts}
\usepackage{defs}
\usepackage{placeins}
\usepackage{booktabs} % for \cmidrule in tables
\usepackage{url}

\newif\ifexternalizing 
	\externalizingfalse

\usepackage{tikz}
\usepackage{pgfplots}
     
\ifexternalizing
     \definecolor{darkgreen}{RGB}{0, 128,0}
     \usepgfplotslibrary{external}
     \pgfplotscreateplotcyclelist{color mod}{%
         blue!30,every mark/.append style={fill=blue!30!black},mark=*\\%
         red!60,every mark/.append style={fill=red!30!black},mark=square*\\%
         brown!60!black,every mark/.append style={fill=brown!80!black},mark=diamond\\%
         black,mark=triangle\\%
         blue,every mark/.append style={fill=blue!80!black},mark=diamond*\\%
         red,densely dashed,every mark/.append style={solid,fill=red!80!black},mark=*\\%
         brown!60!black,densely dashed,every mark/.append style={
         solid,fill=brown!80!black},mark=square*\\%
         black,densely dashed,every mark/.append style={solid,fill=gray},mark=otimes*\\%
         blue,densely dashed,mark=star,every mark/.append style=solid\\%
         red,densely dashed,every mark/.append style={solid,fill=red!80!black},mark=diamond*\\%
     }
\else
\fi

\definecolor{greenyellow}   {cmyk}{0.15, 0   , 0.69, 0   }
\definecolor{yellow}        {cmyk}{0   , 0   , 1   , 0   }
\definecolor{goldenrod}     {cmyk}{0   , 0.10, 0.84, 0   }
\definecolor{dandelion}     {cmyk}{0   , 0.29, 0.84, 0   }
\definecolor{apricot}       {cmyk}{0   , 0.32, 0.52, 0   }
\definecolor{peach}         {cmyk}{0   , 0.50, 0.70, 0   }
\definecolor{melon}         {cmyk}{0   , 0.46, 0.50, 0   }
\definecolor{yelloworange}  {cmyk}{0   , 0.42, 1   , 0   }
\definecolor{orange}        {cmyk}{0   , 0.61, 0.87, 0   }
\definecolor{burntorange}   {cmyk}{0   , 0.51, 1   , 0   }
\definecolor{bittersweet}   {cmyk}{0   , 0.75, 1   , 0.24}
\definecolor{redorange}     {cmyk}{0   , 0.77, 0.87, 0   }
\definecolor{mahogany}      {cmyk}{0   , 0.85, 0.87, 0.35}
\definecolor{maroon}        {cmyk}{0   , 0.87, 0.68, 0.32}
\definecolor{brickred}      {cmyk}{0   , 0.89, 0.94, 0.28}
\definecolor{red}           {cmyk}{0   , 1   , 1   , 0   }
\definecolor{orangered}     {cmyk}{0   , 1   , 0.50, 0   }
\definecolor{rubinered}     {cmyk}{0   , 1   , 0.13, 0   }
\definecolor{wildstrawberry}{cmyk}{0   , 0.96, 0.39, 0   }
\definecolor{salmon}        {cmyk}{0   , 0.53, 0.38, 0   }
\definecolor{carnationpink} {cmyk}{0   , 0.63, 0   , 0   }
\definecolor{magenta}       {cmyk}{0   , 1   , 0   , 0   }
\definecolor{violetred}     {cmyk}{0   , 0.81, 0   , 0   }
\definecolor{rhodamine}     {cmyk}{0   , 0.82, 0   , 0   }
\definecolor{mulberry}      {cmyk}{0.34, 0.90, 0   , 0.02}
\definecolor{redviolet}     {cmyk}{0.07, 0.90, 0   , 0.34}
\definecolor{fuchsia}       {cmyk}{0.47, 0.91, 0   , 0.08}
\definecolor{lavender}      {cmyk}{0   , 0.48, 0   , 0   }
\definecolor{thistle}       {cmyk}{0.12, 0.59, 0   , 0   }
\definecolor{orchid}        {cmyk}{0.32, 0.64, 0   , 0   }
\definecolor{darkorchid}    {cmyk}{0.40, 0.80, 0.20, 0   }
\definecolor{purple}        {cmyk}{0.45, 0.86, 0   , 0   }
\definecolor{plum}          {cmyk}{0.50, 1   , 0   , 0   }
\definecolor{violet}        {cmyk}{0.79, 0.88, 0   , 0   }
\definecolor{royalpurple}   {cmyk}{0.75, 0.90, 0   , 0   }
\definecolor{blueviolet}    {cmyk}{0.86, 0.91, 0   , 0.04}
\definecolor{periwinkle}    {cmyk}{0.57, 0.55, 0   , 0   }
\definecolor{cadetblue}     {cmyk}{0.62, 0.57, 0.23, 0   }
\definecolor{cornflowerblue}{cmyk}{0.65, 0.13, 0   , 0   }
\definecolor{midnightblue}  {cmyk}{0.98, 0.13, 0   , 0.43}
\definecolor{navyblue}      {cmyk}{0.94, 0.54, 0   , 0   }
\definecolor{royalblue}     {cmyk}{1   , 0.50, 0   , 0   }
\definecolor{blue}          {cmyk}{1   , 1   , 0   , 0   }
\definecolor{cerulean}      {cmyk}{0.94, 0.11, 0   , 0   }
\definecolor{cyan}          {cmyk}{1   , 0   , 0   , 0   }
\definecolor{processblue}   {cmyk}{0.96, 0   , 0   , 0   }
\definecolor{skyblue}       {cmyk}{0.62, 0   , 0.12, 0   }
\definecolor{turquoise}     {cmyk}{0.85, 0   , 0.20, 0   }
\definecolor{tealblue}      {cmyk}{0.86, 0   , 0.34, 0.02}
\definecolor{aquamarine}    {cmyk}{0.82, 0   , 0.30, 0   }
\definecolor{bluegreen}     {cmyk}{0.85, 0   , 0.33, 0   }
\definecolor{emerald}       {cmyk}{1   , 0   , 0.50, 0   }
\definecolor{junglegreen}   {cmyk}{0.99, 0   , 0.52, 0   }
\definecolor{seagreen}      {cmyk}{0.69, 0   , 0.50, 0   }
\definecolor{green}         {cmyk}{1   , 0   , 1   , 0   }
\definecolor{forestgreen}   {cmyk}{0.91, 0   , 0.88, 0.12}
\definecolor{pinegreen}     {cmyk}{0.92, 0   , 0.59, 0.25}
\definecolor{limegreen}     {cmyk}{0.50, 0   , 1   , 0   }
\definecolor{yellowgreen}   {cmyk}{0.44, 0   , 0.74, 0   }
\definecolor{springgreen}   {cmyk}{0.26, 0   , 0.76, 0   }
\definecolor{olivegreen}    {cmyk}{0.64, 0   , 0.95, 0.40}
\definecolor{rawsienna}     {cmyk}{0   , 0.72, 1   , 0.45}
\definecolor{sepia}         {cmyk}{0   , 0.83, 1   , 0.70}
\definecolor{brown}         {cmyk}{0   , 0.81, 1   , 0.60}
\definecolor{tan}           {cmyk}{0.14, 0.42, 0.56, 0   }
\definecolor{gray}          {cmyk}{0   , 0   , 0   , 0.50}
\definecolor{black}         {cmyk}{0   , 0   , 0   , 1   }
\definecolor{white}         {cmyk}{0   , 0   , 0   , 0   } 

 \usepackage[plainpages=false,pdfpagelabels,bookmarksnumbered,%
        colorlinks=true,%
        linkcolor=sepia,%
        citecolor=sepia,%
        filecolor=maroon,%
        pagecolor=red,%
        urlcolor=sepia]{hyperref} 

\newcommand{\externaltikz}[2]{
\ifexternalizing
     \tikzsetnextfilename{#1} #2

  \pgfplotsset{compat=newest}
  \pgfplotsset{plot coordinates/math parser=false}
  
\else
     \includegraphics{Externals/#1}
\fi
}

    \newlength\figureheight
    \newlength\figurewidth
\newtheorem{theorem}{Theorem}[section]
\newtheorem{definition}[theorem]{Definition}
\newtheorem{remark}[theorem]{Remark}

\newcommand{\ceil}[1]{\ensuremath{\left\lceil #1 \right\rceil}}

\newcommand{\secref}[1]{Section~\ref{#1}}
\newcommand{\thmref}[1]{Theorem~\ref{#1}}

\newcommand{\figref}[1]{Figure~\ref{#1}}
\newcommand{\tabref}[1]{Table~\ref{#1}}

\newcommand{\R}{\mathbb{R}}

\newcommand{\U}{\ensuremath{\mathbf{u}}}
\newcommand{\RD}[2]{\ensuremath{\mathcal{R}_{#1}^{#2}}}
\newcommand{\RQ}[1]{\RD{#1}{\mathcal{Q}}}

\newcommand{\intQ}[1]{\ensuremath{\left<#1\right>}}
\newcommand{\nmom}{\ensuremath{N}}

\newcommand{\nqmu}{\ensuremath{{N_\mathcal{Q}}}}
\newcommand{\nqs}{\ensuremath{{Q}}}
\newcommand{\order}{\ensuremath{k}}
\newcommand{\orderS}{\order_{\source}}
\newcommand{\timelvl}{\ensuremath{n}}
\newcommand{\ngrid}{\ensuremath{{J}}}

\newcommand{\Vintp}[1]{\Vint{#1}_+}
\newcommand{\Vintm}[1]{\Vint{#1}_-}

\newcommand{\basisComp}{\ensuremath{m}}
\newcommand{\basis}{\ensuremath{ \mathbf{\basisComp}}}

\newcommand{\PN}{\ensuremath{\text{P}_{\nmom}}}
\newcommand{\MN}{\ensuremath{\text{M}_{\nmom}}}

\newcommand{\nQ}{n_{\cQ}}

\DeclareMathOperator*{\argmin}{argmin}

\newcommand{\source}{S}

\def\xL{x_{\text{L}}}
\def\xR{x_{\text{R}}}
\def\psiL{\psi_{\text{L}}}
\def\psiR{\psi_{\text{R}}}
\def\tf{t_{\rm f}}
\def\psiInit{\psi_{t = 0}}
\def\etad{\eta_{\ast}}

\def\alphahat{\hat{\bsalpha}}
\def\alphabar{\overline{\bsalpha}}

\def\uIso{\bu_{\text{iso}}}
\def\collision{\mathcal{C}}

\def\xL{x_{\text{L}}}
\def\xR{x_{\text{R}}}
\def\psiL{\psi_{\text{L}}}
\def\psiR{\psi_{\text{R}}}
\def\tf{t_{\rm f}}
\def\psiInit{\psi_{t = 0}}
\def\etad{\eta_{\ast}}

\def\ansatz{\hat \psi}

\newcommand{\psibar}{\overline{\psi}}
\newcommand{\ubar}{\overline{\bu}}

\numberwithin{equation}{section}
\title{A realizability-preserving high-order kinetic scheme using
WENO reconstruction for entropy-based moment closures of linear kinetic
equations in slab geometry}
\author[fs]{Florian Schneider}
\author[ga]{Graham Alldredge}
\author[jk]{Jochen Kall}
\address[fs]{Fachbereich Mathematik, TU Kaiserslautern, Erwin-Schr\"odinger-Str., 67663 Kaiserslautern, Germany, {\tt schneider@mathematik.uni-kl.de}}
\address[ga]{Department of Mathematics, RWTH Aachen University, Schinkelstr. 2,52062 Aachen, Germany, {\tt alldredge@mathcces.rwth-aachen.de}}
\address[jk]{Fachbereich Mathematik, TU Kaiserslautern, Erwin-Schr\"odinger-Str., 67663 Kaiserslautern, Germany, {\tt kall@mathematik.uni-kl.de}}

\date{}

\begin{document}
\begin{abstract}
We develop a high-order kinetic scheme for entropy-based moment models of
a one-dimensional linear kinetic equation in slab geometry.
High-order spatial reconstructions are achieved using the weighted
essentially non-oscillatory (WENO) method, and for time integration we use
multi-step Runge-Kutta methods which are strong stability preserving and
whose stages and steps can be written as convex combinations of forward
Euler steps.
We show that the moment vectors stay in the realizable set using these time
integrators along with a maximum principle-based kinetic-level limiter,
which simultaneously dampens spurious oscillations in the numerical
solutions.
We present numerical results both on a manufactured solution, where we
perform convergence tests showing our scheme converges of the expected
order up to the numerical noise from the numerical optimization,
as well as on two standard benchmark problems, where we show some of the
advantages of high-order solutions and the role of the key parameter in the
limiter.
\end{abstract}
\maketitle

\noindent

\input{Sections/Introduction}
\input{Sections/Macroscopic}
\input{Sections/KineticScheme}
\input{Sections/limiting}
\input{Sections/Results}
\input{Sections/Conclusions}

%%%%%%%%%%%%
% Bibliography
%%%%%%%%%%%%%%
\bibliographystyle{plain}
\bibliography{RadLit,Qiqqa2BibTexExport}

\end{document}

%% file: Sections/Introduction.tex
%!TEX root = ../MMEOptimization.tex

%%%%%%%%%%%%%%%%%%%%%%
\section{Introduction}
%%%%%%%%%%%%%%%%%%%%%%
In recent years many approaches have been considered for the solution of
time-dependent linear kinetic transport equations, which arise for example
in electron radiation therapy or radiative heat transfer problems.
Many of the most popular methods are moment methods, also known as moment
closures because they are distinguished by how they close the truncated
system of exact moment equations.
Moments are defined through angular averages against basis functions to
produce spectral approximations in the angle variable.
A typical family of moment models are the so-called P${}_\nmom$-methods
\cite{Lewis-Miller-1984,Gel61} which are pure spectral methods.
However, many high-order moment methods, including P${}_\nmom$, do not take
into account that the original kinetic density to be approximated must be
nonnegative.
The moment vectors produced by such models are therefore often not realizable,
that is, there is no associated nonnegative kinetic density consistent with
the moment vector, and thus the solutions can contain obviously non-physical
artifacts such as negative local particle densities \cite{Bru02}.

The family of minimum-entropy models, colloquially known as $\MN$ models
or entropy-based moment closures, solve this problem (for certain physically
relevant entropies) by specifying the closure using a nonnegative density
reconstructed from the moments.
The $\MN$ models are the only models which additionally are hyperbolic and
dissipate entropy \cite{Lev96}.
The cost of all these properties is that the reconstruction of this density
involves solving an optimization problem at every point on the space-time
mesh.
These reconstructions, however, can be parallelized, and so the recent emphasis
on algorithms that can take advantage of massively parallel computing
environments has led to renewed interest in the computation of $\MN$ solutions
both for linear and nonlinear kinetic equations
\cite{DubFeu99,Hauck2010,LamGroth2011,AllHau12,GarrettHauck2013,
McDonaldGroth2013}.
Despite the parallelizability of the cost of the numerical optimization, the
gain in efficiency that would come from a higher-order space-time
discretization will still be necessary for a practical $\MN$ implementation.

The key challenge for high-order methods for entropy-based moment closures is
that the numerical solutions leave the set of realizable moments
\cite{Olbrant2012}, outside of which the defining optimization problem has no
solution.
Discontinuous-Galerkin methods can handle this problem using a realizability
limiter directly on the moment vectors themselves
\cite{Zhang2010,Olbrant2012,AlldredgeSchneider2014}, but at this level
realizability conditions are in general quite complicated
and also not well-understood for two- or three-dimensional problems for
moment models of order higher than two.
Realizability limiting for kinetic schemes, however, is much easier because
at the level of the kinetic density, realizability corresponds simply to
nonnegativity.
Furthermore, this same limiter can be strengthened to also enforce a local
maximum principle, thereby dampening artificial oscillations in numerical
solutions.

Thus in this work we derive a high-order (in space and time) kinetic scheme
for $\MN$ models with moments of (in principle) arbitrary order.
We start in \secref{sec:momapprox} by introducing the linear kinetic equation
we will consider, its entropy-based moment closure, and reviewing the concept
of realizability.
Continuing in \secref{sec:RPScheme} we introduce the concept of a kinetic
scheme for moment equations and then give our numerical techniques for the 
discretization of each of the independent variables: angle, space, and time.
The issue of realizability preservation and the necessary limiters are
discussed in \secref{sec:limiting}, finishing the full description of our
scheme.
The results from our numerical simulations are presented in
\secref{sec:Results}, including a convergence study using a manufactured
solution and solutions for two benchmark problems.
Finally we draw conclusions and discuss the next steps for future work in
\secref{sec:Conclusions}.

%% file: Sections/Macroscopic.tex
%!TEX root = ../MMEOptimizationKRM.tex

\section{A linear kinetic equation and moment closures}
\label{sec:momapprox}

We begin with the linear kinetic equation we will use to test our algorithm
and a brief introduction to entropy-based moment closures which closely
follows \cite{AlldredgeSchneider2014}.
More background can be found for example in \cite{Lewis-Miller-1984,
Lev96,Hauck2010} and references therein.

\subsection{A linear kinetic equation}

We consider the following one-dimensional linear kinetic equation for the 
kinetic density $\psi = \psi(t, x, \mu) \ge 0$ in slab geometry,
for time $t > 0$, spatial coordinate
$x \in X = (\xL, \xR) \subseteq \R$, and angle variable
$\mu \in [-1, 1]$:
\begin{align}
\label{eq:FokkerPlanck1D}
\partial_t\psi + \mu \partial_x  \psi + \sig{a} \psi = 
\sig{s}\collision(\psi) + \source,
\end{align}
where $\sig{a}$ and $\sig{s}$ are the nonnegative absorption and scattering
interaction  coefficients, and $\source$ a source.
The operator $\collision$ is a collision operator, which in this paper we
assume to be linear and have the form
\begin{equation}
 \collision(\psi) = \int_{-1}^1 T(\mu, \mu^\prime)
  \psi(t, x, \mu^\prime)~d\mu^\prime 
  - \int_{-1}^1 T(\mu^\prime, \mu) \psi(t, x, \mu)~d\mu^\prime.
\label{eq:collisionOperator}
\end{equation}
We assume that the kernel $T$ is strictly positive and normalized to 
$\int_{-1}^1 T(\mu^\prime, \mu) d\mu^\prime~\equiv~1$.  A typical example is 
isotropic scattering, where $T(\mu, \mu^\prime) \equiv 1/2$.

Equation \eqref{eq:FokkerPlanck1D} is supplemented by initial and boundary 
conditions:
\begin{subequations}
\label{eq:bc-ic}
\begin{align}
 \psi(t, \xL, \mu) &= \psiL(t, \mu) \,, & t &\geq 0 \,, & \mu &> 0\,,
  \label{eq:bcL} \\
 \psi(t, \xR, \mu) &= \psiR(t, \mu) \,, & t &\geq 0 \,, & \mu &< 0\,,
  \label{eq:bcR} \\
 \psi(0, x, \mu) &= \psiInit(x, \mu)  \,, & x &\in 
  (\xL,\xR) \,, & \mu &\in [-1,1] \,,
\end{align}
\end{subequations}
where $\psiL$, $\psiR$, and $\psiInit$ are given.

\subsection{Moment equations and entropy-based closures}
Moment equations are an angular discretization for
\eqref{eq:FokkerPlanck1D}, where the moments themselves are defined by
angular averages against a set of basis functions.  
We use the following notation for angular integrals:
$$
 \Vint{\phi} = \int_{-1}^1 \phi(\mu) d\mu
$$
for any integrable function $\phi = \phi(\mu)$;
and therefore if we collect the basis functions into a vector $\basis = 
\basis(\mu) = (\basisComp_0(\mu), \basisComp_1(\mu), \ldots, 
\basisComp_{\nmom}(\mu))^T$, the moments of a kinetic density
$\phi$ are given by $\U = \Vint{\basis \phi}$.
In this paper we consider the monomial moments $\basisComp_i(\mu)
= \mu^i$, though all results can be extended to other bases, including, for
example, partial \cite{DubKla02,DubFraKlaTho03} or mixed moments \cite{Frank07,SchneiderAlldredge14}.

The closed system of moment equations is a system of partial differential
equations of the form
\begin{equation}
 \partial_t \U + \partial_x \bff(\U) + \sig{a}\U = \sig{s} \br(\U) + 
  \Vint{\basis \source},
\label{eq:moment-closure}
\end{equation}
where the moment vector $\U(t, x)$ approximates $\vint{\basis \psi}$ for the 
kinetic density $\psi$ satisfying \eqref{eq:FokkerPlanck1D}.
In an entropy-based closure (commonly referred to as the $\MN$ model or the
Levermore closure after he exposed their
general structure in \cite{Lev96}), the functions $\bff$ and $\br$
have the form
$$
 \bff(\U) := \Vint{\mu \basis \ansatz_{\U}} \quand
 \br(\U)  := \Vint{\basis \collision(\ansatz_{\U})}.
$$
Here $\ansatz_{\U}$ is an \textit{ansatz} density reconstructed from the
moments  $\U$ by solving the constrained optimization problem:
\begin{equation}
 \ansatz_{\U} = \argmin\limits_\phi \left\{\Vint{\eta(\phi)}
 : \Vint{\basis \phi} = \U \right\},
\label{eq:primal}
\end{equation}
where the kinetic entropy density $\eta$ is strictly convex and
the minimum is simply taken over functions $\phi = \phi(\mu)$ such that 
$\Vint{\eta(\phi)}$ and $\vint{\basis \phi}$ are well defined.
This problem is typically solved through its strictly convex, unconstrained, finite-dimensional dual,
\begin{equation}
 \alphahat(\U) := \argmin_{\bsalpha \in \R^{\nmom + 1}} \Vint{\eta_*(\basis^T 
  \bsalpha)} - \U^T \bsalpha,
\label{eq:dual}
\end{equation}
where $\eta_*$ is the Legendre dual of $\eta$.
The first-order necessary conditions for $\alphahat(\U)$ 
show that the solution to \eqref{eq:primal} has the form
\begin{equation}
 \ansatz_{\U} = \etad' \left(\basis^T \alphahat(\U) \right)
\label{eq:psiME}
\end{equation}
where $\etad'$ is the derivative of $\eta_*$.

The kinetic entropy density $\eta$ can be chosen according to the 
physics being modeled.

While in a linear setting such as ours, indeed any convex entropy $\eta$ is 
dissipated by \eqref{eq:FokkerPlanck1D}. As in \cite{Hauck2010} we focus on the Maxwell-Boltzmann entropy,
\begin{align*}
%\label{eq:EntropyM}
 \eta(z) = z \log(z) - z,
\end{align*}
because not only is it physically relevant for a wide variety of problems,
but in particular it gives a positive ansatz $\ansatz_{\U}$, since
$\etad(y) = \etad'(y) = \exp(y)$, and thus
\begin{align}
\label{eq:mb-ansatz}
 \ansatz_{\U} = \exp \left( \basis^T \alphahat(\U) \right).
\end{align}

Another standard choice for the entropy is $\eta(z) = \frac12 z^2$, which
yields the well-known $\PN$ equations \cite{Lewis-Miller-1984,
Hauck2010}.
The $\PN$ equations are linear, and since $\etad'(y) = y$, the optimization
problem \eqref{eq:dual} can be solved by hand.
However the resulting ansatz is simply a linear combination of the basis
polynomials and thus is not necessarily nonnegative.%
\footnote{
For this reason, while the rest of the scheme we give below can be applied
to the $\PN$ equations, the positivity-preserving techniques we use
(see \secref{sec:limiting}) would be unnecessary.
However, the scheme should apply equally well for any
entropy-based closure with a positive ansatz.
}

The incorporation of the boundary conditions \eqref{eq:bc-ic} is neither
obvious nor trivial
and is still an open problem \cite{Pomraning-1964,Larsen-Pomraning-1993-1, 
Larsen-Pomraning-1993-2,Struchtrup-2000}.
This is not a focus of our work here, and therefore we only use a
simple approach from previous work on entropy-based moment closures
\cite{Hauck2010,AllHau12}, which we discuss below in \secref{sec:WENO}.

\subsection{Moment realizability}

Since the underlying kinetic density we are trying to approximate is
nonnegative, a 
moment vector only makes sense physically if it can be associated with a 
nonnegative density. In this case the moment vector is called 
\textit{realizable}.
Additionally, since the entropy ansatz has the form \eqref{eq:mb-ansatz},
the optimization problem \eqref{eq:primal} only has a 
solution if the moment vector lies in the ansatz moment space
$$
 \cA := \left\{ \Vint{\basis \exp \left( \basis^T \bsalpha \right)} 
  : \bsalpha \in \bbR^{N + 1}  \right\}.
$$
In our case, where the domain of angular integration is bounded, the ansatz 
moment space $\cA$ is exactly equal to the set of realizable moment vectors 
\cite{Jun00}.
Therefore we can focus simply on realizable moments:

\begin{definition}
The \emph{realizable set} \RD{\basis}{} is 
$$
\RD{\basis}{} = \left\{\U~:~\exists \phi(\mu)\ge 0,\, \Vint{\phi} > 0,
 \text{ such that } \U = \Vint{\basis\phi} \right\}.
$$
Any $\phi$ such that $\U = \vint{\basis \phi}$ is called a \emph{representing 
density}.
\end{definition}

The realizable set is a convex cone.
In the monomial basis, a moment vector is realizable if 
and only if its corresponding Hankel matrices are positive definite 
\cite{Shohat-Tamarkin-1943,CurFial91}.

In general, angular integrals cannot be computed analytically.
We define a quadrature for functions $\phi:[-1,1]\to\R$ by nodes
$\{\mu_i\}_{i=1}^{\nqmu}$ and weights $\{w_i\}_{i=1}^{\nqmu}$ such that
\begin{align*}
\sum\limits_{i=1}^{\nqmu} w_i \phi(\mu_i) \approx \Vint{\phi}
\end{align*}
Then the numerically realizable set is \cite{ahot2013}
$$
\RQ{\basis} = \left\{\U~:~\exists f_i > 0 \text{ s.t. } \U = \sum_{i = 
1}^{\nqmu}w_i \basis(\mu_i) f_i \right\}\subset\RD{\basis}{} %\label{eq:QRealizable}
$$
Indeed, when replacing the integrals in the optimization problem \eqref{eq:primal} 
with quadrature, a minimizer can only exist when $\U \in \RQ{\basis}$.
Below we often abuse notation and write $\Vint{\phi}$ when in 
implementation we mean its approximation by quadrature.
We also liberally use the term \emph{realizable} either to mean realizability
with respect to $\RD{\basis}{}$ or with respect to $\RQ{\basis}$, where
the specific meaning depends on whether exact integrals or those
approximated by quadrature are meant in the context.

%% file: Sections/KineticScheme.tex
%!TEX root = ../MMEOptimizationKRM.tex
% !TeX spellcheck = en_US
%%%%%%%%%%%%%%%%%%%%%%
\section{A high-order kinetic scheme}
%%%%%%%%%%%%%%%%%%%%%%
\label{sec:RPScheme}
A kinetic scheme for \eqref{eq:moment-closure} can be thought of as first
defining a spatial discretization for the underlying kinetic equation
\eqref{eq:FokkerPlanck1D} and subsequently performing the angular
discretization with the moment closure.
We divide the spatial domain $(\xL, \xR)$ into a (for simplicity) uniform 
grid of $\ngrid$ cells $I_j = (x_{j - 
1/2}, x_{j + 1/2})$, where the cell edges are given by $x_{j \pm 1/2} = x_j \pm 
\dx / 2$ for cell centers $x_j = \xL + (j - 1 / 2)\dx$, and
$\dx = (\xR - \xL) /  \ngrid$. 
For \eqref{eq:FokkerPlanck1D} we define a finite-volume scheme for the cell
means
\begin{equation}
\psibar_j(t, \mu) \simeq \frac{1}{\dx} \int_{I_j} \psi(t, x, \mu) dx,
\label{eq:psi-cell-avg}
\end{equation}
which with the Godunov (or `upwind') numerical flux gives:
\begin{align*}
 \partial_t \psibar_j &+ \max(\mu, 0) \cfrac{\psi_{j + 1 / 2}^-
  - \psi_{j - 1 / 2}^-}{\Delta x} \\
 &+ \min(\mu, 0) \cfrac{\psi_{j + 1 / 2}^+
  - \psi_{j - 1 / 2}^+}{\Delta x} + \overline{\sigma_a\psi}_j
  = \frac12 \overline{\sig{s} \collision(\psi)}_j + \overline{\source}_j
\end{align*}
where $\psi_{j \pm 1 / 2}^-$ and $\psi_{j \pm 1 / 2}^+$ in the flux terms
denote the values of
the approximate solution at the cell edges $x_{j \pm 1/2}$ from the
left and right, respectively, and we generally use
the bar with subsequent subscript $j$, i.e. $\overline{\, \cdot \,}_j$,
to indicate a cell average over the $j$-th cell as in
\eqref{eq:psi-cell-avg}.

To obtain a high-order scheme in space one only has to give a high-order
reconstruction of the point-values of $\psi$, the distribution 
underlying the cell means, not only at the edge values for the flux terms,
but also throughout the cells when $\sig{a}$ or $\sig{s}$ depends on
$x$.
We use the popular weighted essentially non-oscillatory (WENO) reconstruction
method \cite{Shu1998}, which gives a polynomial reconstruction
of $\psi$ from the cell averages of the $j$-th cell and its neighbors.

Now we perform the moment closure by replacing $\psi$ with the entropy ansatz $\ansatz$ in \eqref{eq:psiME},
multiplying through by the angular basis functions, and integrating out the
angle:

\begin{align}
 \partial_t \ubar_j + \Vintp{\mu \basis
  \cfrac{\ansatz_{j + 1 / 2}^-
  - \ansatz_{j - 1 / 2}^-}{\dx}} 
  + \Vintm{\mu \basis \cfrac{\ansatz_{j + 1 / 2}^+
  - \ansatz_{j - 1 / 2}^+}{\dx}}
  + \overline{\sigma_a \bu}_j
  = \cfrac{1}{2}\overline{\sig{s} \br(\bu)}_j + \overline{\bs}_j
\label{eq:semi-discrete-mom-eqn}
\end{align} 
where
$$
 \Vintp{\phi} = \int_0^1 \phi(\mu) d\mu \quand
 \Vintm{\phi} = \int_{-1}^0 \phi(\mu) d\mu,
$$
and $\bs = \vint{\basis \source}$.
Here $\ansatz_{j \pm 1 / 2}^+$ denote the evaluations at the cell edges of
the WENO reconstructions made using the entropy ans\"atze of the neighboring
cells evaluated from the right, and respectively for $\ansatz_{j \pm 1 / 2}^-$
evaluated from the left.
The first step in the scheme, then, is to compute the ansatz for
each cell.

\subsection{Numerical optimization for angular reconstruction}
\label{sec:Optimization}

In order to compute $\ansatz_j$ at the cell means,
we first compute the multipliers
$\alphahat(\ubar_j)$ by solving the dual problem \eqref{eq:dual}.
The gradient and Hessian of the objective function are
\begin{equation}
 \bg(\bsalpha) = \Vint{\basis \exp(\basis^T \bsalpha)} - \ubar_j
  \quand
 \bH(\bsalpha) = \Vint{\basis \basis^T \exp(\basis^T \bsalpha)}.
\label{eq:grad-hess}
\end{equation}

We use the numerical optimization techniques proposed in \cite{ahot2013}.

We assume that the moment vector $\ubar_j$ has been scaled so that its
zeroth component is one.
The optimizer stops at the first iterate $\bsalpha$ which satisfies
\begin{subequations}
\label{eq:stop}
\begin{align}
 \| \bg(\bsalpha) \|
  & < \frac{\tau}{1 + \| \ubar_j \| + \tau} =: \tau', \text{ and}
  \label{eq:stop-grad} \\
 1 - \veps &<  \exp(-\| \bd(\bsalpha) \|_1 - |\log (u_0(\bsalpha))|) 
  \label{eq:stop-gamma}
\end{align}
\end{subequations}
where $\| \cdot \|$ is the Euclidean and $\| \cdot \|_1$ the $1-$norm in
$\R^{\nmom + 1}$, $\tau$ and $\veps$ user-specified tolerances,
and $u_0(\bsalpha) = \vint{\exp(\basis^T \bsalpha}$ is the zero-th order
moment associated with the multiplier vector $\bsalpha$.
This criterion is similar to the one in \cite{ahot2013} but modified for the following
reasons.

Unlike the algorithm there, we modify the zero-th component of
the final multipliers so that the zero-th moment (which gives the local
density) is exactly matched.
That is, while the optimizer stops at the first $\bsalpha$ which satisfies
\eqref{eq:stop}, the multiplier vector it returns (to define $\psibar_j$
for the kinetic reconstructions in \secref{sec:WENO}) is
\begin{equation}
\label{eq:alphabar}
\alphabar := \begin{pmatrix} \alpha_0 - \log(u_0(\bsalpha)),
 \alpha_1, \dots, \alpha_{\nmom}
 \end{pmatrix}^T.
\end{equation}
Since $\basisComp_0 \equiv 1$, this gives $\vint{\basis
\exp(\basis^T \alphabar)} = \vint{\basis \exp(\basis^T \bsalpha)}
/ u_0(\bsalpha)$, which ensures $\vint{\exp(\basis^T \alphabar)} = 1$.
The form of $\tau'$ on the right-hand side of
\eqref{eq:stop-grad} ensures that $\|\bg(\alphabar)\|$ is bounded by $\tau$:
\begin{align}
\label{eq:GradientNorm}\| \bg(\alphabar) \| 
 & = \| \Vint{\basis \exp(\basis^T \alphabar)} - \ubar_j \| \\
 & \le \frac{1}{u_0(\bsalpha)} \left\| \Vint{\basis
  \exp(\basis^T \bsalpha)} - \ubar_j \right\|\nonumber
   + \left\| \frac{1}{u_0(\bsalpha)}\ubar_j - \ubar_j \right\| \\
 & \le \frac{1}{1 - \tau'} \tau'
   + \frac{\tau'}{1 - \tau'} \| \ubar_j \|
  = \tau, \nonumber
\end{align}
where we have used that \eqref{eq:stop-grad} implies $|u_0(\bsalpha) - 1|
< \tau'$ and consequently $1 / u_0(\bsalpha) < 1 / (1 - \tau')$.

The second condition \eqref{eq:stop-gamma} in the stopping criterion
enforces an approximate lower bound on the ratio%
\footnote{This is also similar to what was done in \cite{ahot2013}, though
there the bound was more naturally given as an upper bound of the inverse of
the ratio we use.
}
$\ansatz_j / \psibar_j$, where $\ansatz_j = \ansatz_{\ubar_j} =
\exp(\basis^T \alphahat(\ubar_j))$
is the entropy ansatz associated with an exact solution of
\eqref{eq:mb-ansatz} and $\psibar_j = \exp(\basis^T \alphabar)$
is the ansatz associated with the Lagrange multiplier vector which satisfies
\eqref{eq:stop}.
While of course we do not know the exact entropy ansatz, we can approximate
the difference between the exact solution $\alphahat(\bu_j)$ and another
multiplier vector $\bsalpha$ from the iterations of the optimizer by the
Newton direction
$$
\bd(\bsalpha) = -\bH(\bsalpha)^{-1} \bg(\bsalpha).
$$
This leads to the approximate bound
\begin{align*}
 \frac{\ansatz_j}{\psibar_j} = \exp(\basis^T(\alphahat(\ubar_j) - \alphabar))
  & =  \exp(\basis^T(\alphahat(\ubar_j) - \bsalpha + \bsalpha - \alphabar)) \\
  & \approx \exp(\basis^T ( \bd(\bsalpha) + \bsalpha - \alphabar)) \\
  & \ge \exp(-\|\basis\|_\infty (\| \bd(\bsalpha) \|_1
    + |\log(u_0(\bsalpha))|)),
\end{align*}
where $\|\bm\|_\infty = \max_{i,\mu} |\basisComp_i(\mu)| = 1$.

Finally, and exactly as in \cite{ahot2013}, we use a isotropic-regularization
technique to return multipliers for nearby moments when the optimizer fails
(for example, by reaching a maximum number of iterations or being
unable to solve for the Newton direction).
Isotropically regularized moments are defined by the convex combination
\begin{align}
\label{eq:RegularizedMoments}
\bv(\U, r) := (1 - r) \U + r u_0 \uIso,
\end{align}

where $\uIso = \frac12 \Vint{\basis}$ is the moment vector of the normalized
isotropic  density $\phi(\mu) \equiv 1/2$.
Then the optimizer moves through a sequence of values $0, r_1,
r_2, \ldots, r_M$, advancing in this sequence only
if the optimizer fails to converge for $\bv(\U, r)$ after $k_{\rm reg}$
iterations for the current value of $r$.
It is assumed that $r_M$ is chosen large enough that the optimizer
will always converge for $\bv(\U, r_M)$ for any realizable $\U$.

\subsection{Spatial WENO reconstruction}
\label{sec:WENO}

\def\cellindex{j}
\def\rightshift{i}
\def\sumindex{m}
\def\linksrechts{\pm}

We use the standard WENO reconstruction method given for example in \cite{Shu1998,toro2009riemann}.
For the unfamiliar reader, in this section we briefly introduce the method,
while a more detailed documentation and demo implementations of the
reconstruction procedures we implemented can be found on our webpage
\cite{AGTMWENOpage}.

Here, as in the previous section, we use $\psibar_j$ to indicate the entropy
ansatz associated with the Lagrange multipliers returned by the optimizer
to approximate the true multipliers $\alphahat(\ubar_j)$.
Time dependence is again suppressed for clarity of exposition.

At $x = x_{j - 1 / 2}$, the cell edge between the $(j - 1)$-th and $j$-th
cells, for each $\mu$ we evaluate a weighted combination of
polynomials of degree $\order - 1$, $p_{jm}(\cdot, \mu) \in \bbP_{\order-1}$
(the space of polynomials up to degree $\order - 1$),
$m = 0, 1, \ldots , \order$, each solving the interpolation problem
\begin{equation}
 \cfrac{1}{\dx} \int_{I_\ell} p_{jm}(x, \mu) dx =  \psibar_\ell(\mu),
  \qquad \ell \in \{j - \order + m, \ldots ,  j + m - 1\}.
\label{eq:interpolation}
\end{equation}
The WENO method then gives weights $\omega^{\linksrechts}_{j - 1 / 2, m}$
to form the weighted averages
$$
p_{j - 1 / 2}^{\linksrechts}(x, \mu) := \sum_{m = 0}^\order
 \omega^{\linksrechts}_{j - 1 / 2, m} p_{jm}(x, \mu),
$$
and finally we approximate the values at the cell interfaces by 
\[
 \psi_{j - 1 / 2}^-(\mu) \simeq p_{j - 1 / 2}^-(x_{j - 1 / 2}, \mu),
  \quand
 \psi_{j - 1 / 2}^+(\mu) \simeq p_{j - 1 / 2}^+(x_{j - 1 / 2}, \mu).
\]
The weights $\omega^{\linksrechts}_{j - 1 / 2, m}$ are non-linear functions
of the cell-averages and reflect the smoothness of each polynomial $p_{jm}$.
They are computed such that for smooth data the approximation order at the 
cell edge is maximized.
This gives an order $2 \order - 1$ approximation at the cell edge, while
the overall order in the interior of the cell is $\order$.

When at least one of the interaction coefficients $\sig{a}$ or
$\sig{s}$ is spatially dependent, we must also specify the reconstruction
inside each cell to compute, for example, the
$\overline{\sig{a} \bu}_j$ term in \eqref{eq:semi-discrete-mom-eqn}.
Here we must make a choice, because both $p_{j - 1/2}^+$ and $p_{j + 1/2}^-$
are order $\order$ reconstructions of the density $\psi$ in the $j$-th cell.%
\footnote{
Some reconstruction methods, such as subcell WENO \cite{cheng2013sub} or minmod \cite{toro2009riemann},
give only one polynomial reconstruction inside each cell, and so for these
methods such a choice would be unnecessary.
}
We denote this polynomial $\psi_j(x,\mu)$ and choose it to be
\begin{equation}
 \psi_j(x, \mu) \simeq \begin{cases}
   p_{j + 1/2}^-(x, \mu) & \mbox{if } \mu > 0, \\
   \frac12 \left( p_{j + 1/2}^-(x, \mu) + p_{j - 1/2}^+(x, \mu) \right)
    & \mbox{if } \mu = 0, \\
   p_{j - 1/2}^+(x, \mu) & \mbox{if } \mu < 0.
  \end{cases}
\label{eq:reconstruction-in-cell}
\end{equation}
This particular reconstruction allows us to derive a
realizability-preserving time step in \thmref{thm:main}.

Finally, it remains to incorporate boundary conditions.
We define `ghost cells' at the cell indices
$j \in \{1 - \order, \ldots , 0, J + 1, \ldots , J + \order \}$, namely those
indices which are used in \eqref{eq:interpolation} but have not yet been
defined.
We assume that we can smoothly extend $\psiL$ and $\psiR$ in $\mu$ to
$[-1,1]$.
We then use the simplest possible approach%
\footnote{In a smooth setting this might reduce the order of accuracy at the
boundary to one.}
and set
\begin{subequations}
\label{eq:numericalboundaryconditions}
\begin{align}
 \psibar_j(t,\mu) &:= 
 \begin{cases}
 \psiL(t,\mu) & \text{ if } j \in \{1 - \order, \ldots , 0 \}\\
 \psiR(t,\mu) & \text{ if } j \in \{J + 1, \ldots , J + \order\}.
 \end{cases}
\end{align}
\end{subequations}

We note, however, that the validity of this approach is not entirely
noncontroversial, but the question of appropriate boundary conditions for
moment models is an open problem
\cite{Pomraning-1964,Larsen-Pomraning-1993-1, 
Larsen-Pomraning-1993-2,Struchtrup-2000}
which we do not explore here. 

In case of periodic boundary conditions, we use the data from the
physical cells $j \in \{J - \order+1, \ldots , J \}$ from the right side of
the domain to fill the ghost cells $j \in \{1 - \order,
\ldots , 0 \}$. The ghost cells on the right side of the domain analogously take
the values from the left side of the physical domain.

% % % % % % % % % % % % % % % % % % % % % % % % % % % % % % % % % 
\subsection{High-order time integration}
\label{sec:HighOrderTime}
If we collect the approximate cell means from each spatial cell into one
long vector $\U_h(t) := (\ubar_1^T(t), \ubar_2^T(t), \ldots ,
\ubar_{\ngrid}^T(t))^T$, then equation \eqref{eq:semi-discrete-mom-eqn} can
be written as
$$
\partial_t \bu_h = L_h(\bu_h)
$$
Since entropy-based moment closures are only defined on the realizable set,
it is important to choose a time integrator for which we can prove that
realizability is preserved.
Therefore we follow \cite{Zhang2010,AllHau12} and use a
strong stability-preserving (SSP) method whose
stages and steps are convex combinations of forward Euler steps.
Since the realizable set is convex, the analysis of a forward Euler step
then suffices to prove realizability preservation of the high-order method.

When possible we use a SSP Runge-Kutta (SSP-RK) method, but such methods only exist up to order four \cite{Ruuth2004,Gottlieb2005}.
For higher orders we use the so-called two-step Runge-Kutta (TSRK) SSP
methods \cite{Ketcheson2011} as well as their generalizations, the
multi-step Runge-Kutta (MSRK) SSP methods \cite{Bresten2013}.
They combine Runge-Kutta schemes with positive weights and
high-order multistep-methods to achieve a total order which is higher than
four while preserving the important SSP property.
If we let $\bu_h^{\timelvl}$ indicate the collection of the numerical
approximations to the cell averages of the solution at the $\timelvl$-th time
instant $t_{\timelvl} = n \dt$,
an $s$-stage TSRK method in the low-storage implementation has
the following form \cite{Ketcheson2011}:

\begin{align*}
 y_\ell &= d_\ell \bu_h^{\timelvl - 1} + \left(1 - d_\ell
  - \sum_{m = 0}^{s} q_{\ell m} \right) \bu_h^{\timelvl}
  + \sum_{m = 0}^{s} q_{\ell m} \left( y_m + \cfrac{\dt}{\rho} L_h(y_m)
  \right), \qquad 0 \leq \ell \leq s, \\
 \bu_h^{\timelvl + 1} &= \zeta \bu_h^{\timelvl - 1}
  + \left(1 - \zeta - \sum_{m = 0}^{s} \eta_m \right) 
  \bu_h^{\timelvl} + \sum_{m = 0}^{s}\eta_m \left(y_m
  + \cfrac{\dt}{\rho} L_h(y_m) \right),
\end{align*}
where the coefficients $d_\ell, q_{\ell m}, \zeta, \eta_m, \rho$
define the scheme.
In this work we only consider explicit schemes, where $q_{\ell m} = 0$ for
$\ell \ge m$.
The positive coefficient $\rho$ is called the radius of absolute
monotonicity and indicates how much we can scale our time step $\dt$ while
fulfilling the CFL condition for forward Euler steps (see
\secref{sec:realizable-cell-means} below).

Such schemes provide reasonably good effective CFL numbers, which are
the ratios $\rho / s$ for each method.

If the forward-Euler method is stable under the timestep restriction
$\dt \leq \dt_0$, the high order scheme with radius of absolute monotonicity
$\rho$ is stable under the time-step restriction $\dt \leq \rho\dt_0$.
A larger effective CFL indicates that, in order to reach a given final time,
the operator $L_h$ will need to be evaluated fewer times.
The explicit Euler method, which serves as the reference for efficiency in
this context, has an effective CFL condition of $1$ ($s = \rho = 1$).

The effective CFL numbers of the integration schemes used here are given in
\tabref{tab:CFL}.

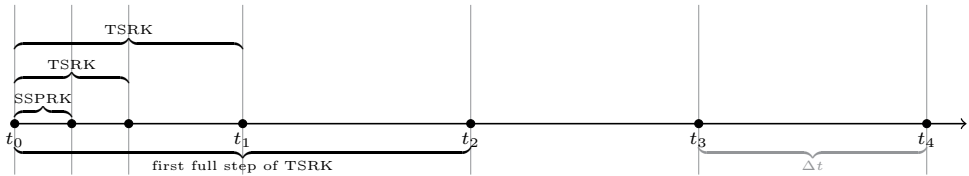
\begin{figure}
  \centering
  \begin{tikzpicture}[scale=0.75]
    \draw [->,semithick] (0,0) -- (16.7,0);

    \foreach \x in {0,1,2,4,8,12,16}
      \draw[gray,very thin] (\x,-0.9) -- (\x, 2.1);
    \foreach \x in {0,1,2,4,8,12,16}
      \filldraw (\x,0) circle (2pt);

    % TODO: this is not the best, if scale is changed above, needs
    % to change here too, better to used named nodes too.
    \draw (4,0) node[below=1ex] {\scriptsize
      $\underbrace{%
        \phantom{\tikz[scale=0.75] \draw (0,0) -- (8,0);}
      }_{\text{first full step of TSRK}}$
    };

    \draw [gray] (14,0) node[below=1ex] {\scriptsize
      $\underbrace{%
        \phantom{\tikz[scale=0.75] \draw (0,0) -- (4,0);}
      }_{\dt}$
    };

    \draw (2,0) node[above=5.5ex] {\scriptsize
      $\overbrace{%
        \phantom{\tikz[scale=0.75] \draw (0,0) -- (4,0);}
      }^{\text{TSRK}}$
    };

    \draw (1,0) node[above=2.5ex] {\scriptsize
      $\overbrace{%
        \phantom{\tikz[scale=0.75] \draw (0,0) -- (2,0);}
      }^{\text{TSRK}}$
    };

    \draw (0.5,0) node[above=-2pt] {\scriptsize
      $\overbrace{%
        \phantom{\tikz[scale=0.75] \draw (0,0) -- (1,0);}
      }^{\text{SSPRK}}$
    };

    \draw (0,0) node[below] {\scriptsize $t_{0}$};
    \draw (4,0) node[below] {\scriptsize $t_{1}$};
    \draw (8,0) node[below] {\scriptsize $t_{2}$};
    \draw (12,0) node[below] {\scriptsize $t_{3}$};
    \draw (16,0) node[below] {\scriptsize $t_{4}$};
  \end{tikzpicture}%
  %}  % end centerline
\caption{One possible startup procedure for SSP TSRK schemes.
The first step from $t_0$ to $t_1$ is subdivided into substeps (here there
are three substeps of sizes $\dt / 4, \dt / 4$, and $\dt / 2$).
A one-step SSP Runge-Kutta scheme is used for the first substep, and
subsequent substeps are taken with the TSRK scheme itself, doubling the
step sizes until reaching $t_1$.
Illustration taken from \cite{Ketcheson2011}.}
\label{fig:TSRKInit}
\end{figure}

Unfortunately, multistep methods are not self-starting, so they need a
predictor for the first time-step.
Since we can only use convex combinations of forward Euler steps for all time
steps in order to prove that realizability is maintained
we must use a lower-order method.
We use the strategy given in \cite{Ketcheson2011}:
First, we predict with a smaller step size $\dt^\star = \dt / 2^q$,
for an integer $q \ge 1$, with the ten-stage, fourth-order explicit
SSPRK$(1,4,10)$ method.
Then we use the corresponding TSRK method and double the step-size after
every iteration until we reach $t = \dt$.
This procedure is shown in \figref{fig:TSRKInit} for $q = 2$.

For the five-step method MSRK$(5,7,12)$, which we use for seventh-order
simulations, we have to initialize four steps.
For these initialization steps we use the two-step method TSRK$(2,7,12)$, whose
initial step we predict using the same method given above.
However the radius of absolute monotonicity $\rho$ for TSRK$(2, 7, 12)$ is
approximately $2.7659$, while for MSRK$(5,7,12)$, the radius of absolute
monotonicity $\rho$ is approximately $3.0886$.
This means the time steps we take after initialization will be longer than
those we can take with the TSRK method during initializating without
violating the realizability-preserving CFL condition.
Therefore we stop increasing the step size in the initialization routine
when we reach $\dt / 2$ (as opposed to $\dt$), and then continue with the
TSRK initialization steps of size $\dt / 2$ until we have computed
$\U_h(t_4) = \U_h(4 \dt)$.
At this point we have all the previous steps we need to compute $\U_h(t_5)$
using the five-stage MSRK$(5,7,12)$ method.

\begin{table}
\centering
\begin{tabular}{c r@{}l r@{.}l}

Order & \multicolumn{2}{c}{Method$(m,\order,s)$} & \multicolumn{2}{c}{Effective CFL} \\

\midrule

$1$ & SSPRK & $(1,1,1)$  & $1$ & $0$  \\
$2$ & SSPRK & $(1,2,20)$ & $0$ & $95$ \\
$3$ & SSPRK & $(1,3,16)$ & $0$ & $75$ \\
$4$ & SSPRK & $(1,4,10)$  & $0$ & $6$\\
$5$ &  TSRK & $(2,5,8)$  & $0$ & $4474$\\
$6$ &  TSRK & $(2,6,12)$ & $0$ & $3653$\\
$7$ &  MSRK & $(5,7,12)$ & $0$ & $3089$\\
 
\end{tabular}
\caption{Methods used here and their effective CFL. The nomenclature
METHOD$(m,\order,s)$ denotes a method with $m$ steps, order $\order$,
and $s$ stages.}
\label{tab:CFL}
\end{table}

%% file: Sections/limiting.tex
%!TEX root = ../MMEOptimizationKRM.tex

\section{Realizability preservation and limiting}
\label{sec:limiting}

The strategy to prove that the moments $\ubar_j$ remain realizable at each
time step and inner stage of the time integrator begins by proving that the kinetic
density reconstructed for the cell means remains nonnegative after a forward
Euler step with a certain time-step restriction.
Since we use time integrators which are a convex combinations of Euler steps
(see \secref{sec:HighOrderTime}) this immediately gives nonnegativity
for all time steps and internal stages of time integration.
This proof, however, requires the assumption that the point-wise values of
the current polynomial reconstruction $\psi_j$ are nonnegative at every
spatial and angular quadrature point.
One therefore introduces a limiter to enforce this nonnegativity.

\subsection{Realizability preservation of the cell means}
\label{sec:realizable-cell-means}

We follow along the lines of the main proof in \cite{AllHau12} and will
provide weaker conditions which follow \cite{Zhang2010,Zhang2011}.

A spatial quadrature rule plays a crucial role here. 
We use Gauss-Lobatto rules which are exact for polynomials of degree
$2\nqs - 3$, where $\nqs$ is the number of quadrature nodes.
These rules are characterized by nodes $y_i$ and weights $\hat{w}_i$ for
$i \in \{ 1, 2, \ldots, \nqs \}$ on the reference interval
$[-\frac12, \frac12]$:
\begin{align}
\label{eq:GaussLobatto}
\int_{-1 / 2}^{1 / 2} f(x)~dx \approx \sum \limits_{i = 1}^{\nqs} \hat{w}_i f(y_i).
\end{align}
We let $x_{ji} := x_j + \dx y_i$ denote the quadrature nodes shifted and
scaled for cell $I_j$, and note that since the Gauss-Lobatto rules include
the endpoints, we have $x_{j1} = x_{j - 1 / 2}$ and $x_{jQ} = x_{j + 1 / 2}$.

We also use the property that, since the collision
kernel $T$ in \eqref{eq:collisionOperator} is nonnegative, there exists a
realizable moment vector
$\U_{\collision}$ such that
\begin{align}
\label{eq:CollisionAssumption}
\intQ{\basis\collision(\psi)} = \U_{\collision} - \U.
\end{align}

\begin{theorem}[Main theorem]
\label{thm:main}
Assume that 
\begin{itemize}
 \item[(i)] for all cells $j \in \{1, 2, \ldots, \ngrid\}$ we have
  $0 \leq \source(t_\timelvl, x)|_{I_j}$, $\sig{a}(t_\timelvl, x)|_{I_j}$,
  and $\sig{s}(t_\timelvl, x)|_{I_j}$ are in $\bbP_{\orderS - 1}(I_j)$;
 \item[(ii)] the cell means $\ubar_j^\timelvl$ at time step $t_{\timelvl}$ are
  realizable;
 \item[(iii)] in each cell the ratio between the exact entropy ansatz and its
  approximation from the optimizer satisfies
  \begin{equation}
   \frac{\ansatz_j(\mu)}{\psibar_j(\mu)} \ge 1 - \veps
  \label{eq:ratio-bound}
  \end{equation}
  for all $\mu\in[-1,1]$ in the angular quadrature set; and
 \item[(iv)] the point-wise values of the polynomial reconstructions
  $\ansatz_j(x, \mu) \in \bbP_{\order - 1}$ at the quadrature nodes of the
  $\nqs$-point Gauss-Lobatto on each cell $I_j$ are nonnegative,
  for $\nqs = \ceil{(\order + \orderS + 1 ) / 2}$.%
  \footnote{Where $\ceil{\cdot}$ is the ceiling function, that is, it
   returns smallest integer bigger than or equal to its argument. Since the Gauss-Lobatto rule is exact for polynomials of degree $2\nqs-3$ this choice guarantees to exactly integrate the occurring polynomials of degree $\left(\order + \orderS-2\right)$.}
\end{itemize}
Let
\begin{align}
\label{eq:sigmatmax}
\sig{t, \max} := \max\limits_{j \in \{ 1, \ldots, \ngrid \} }
 \max\limits_{i \in \{ 1, \ldots, \nqs \} } \sig{s}(t_\timelvl, x_{ji})
 + \sig{a}(t_\timelvl, x_{ji}).
\end{align}

Then under the CFL condition
\begin{align}
\label{eq:CFL}
 \dt \leq (1 - \veps) \min \left( \cfrac{1}{\sig{t, \max}},
  \cfrac{\dx \hat{w}_Q}{1
  + \dx \hat{w}_Q \sig{t, \max}} \right),
\end{align}
the cell means $\ubar_j^{\timelvl+1}$ after one forward Euler step are realizable.
\end{theorem}
\begin{proof}
For simplicity we will neglect the time-index for quantities from $t = t_\timelvl$
and use it only for time-level $t = t_{\timelvl+1}$.
An Euler step is given by
\begin{subequations}
\label{eq:moment-euler-step}
\begin{align}
 \ubar_j^{\timelvl + 1} = \ubar_j &- \frac{\dt}{\dx} \left(
  \Vintp{\mu \basis \left( \psi_{j + 1 / 2}^-
   - \psi_{j - 1 / 2}^- \right)}
  + \Vintm{\mu \basis \left(\psi_{j + 1 / 2}^+
   - \psi_{j - 1 / 2}^+ \right)} \right) \\ 
  & \dt \left(- \overline{\sig{a} \bu}_j
   + \overline{\sig{s} \br(\U)}_j + \overline{\bs}_j \right)
\end{align}
\end{subequations}
and consequently we have $\ubar_j^{\timelvl + 1} = \vint{\basis \phi_j}$ for
$\phi_j$ given by
\begin{subequations}
\label{eq:kinetic-euler-step}
\begin{align}
 \phi_j &= \ansatz_j - \cfrac{\dt}{\dx} \left( \max(\mu, 0) \left(
  \psi_{j + 1 / 2}^- - \psi_{j - 1 / 2}^- \right) + \min(\mu, 0)
  \left( \psi_{j + 1 /2 }^+ - \psi_{j - 1 / 2}^+ \right) \right) \\
 & \quad + \dt \left( -\overline{\sig{t} \psi}_j + \overline{\sig{s}
  \psi_{\collision}}_j + \overline{\source}_j \right)
\end{align}
\end{subequations}
where the total interaction coefficient is defined as $\sig{t} := \sig{a} +
\sig{s}$, and $\psi_{\collision}\geq 0$ is the
entropy ansatz corresponding to the realizable part $\U_{\collision}$ of
the collision operator (see \eqref{eq:CollisionAssumption}).
Note that $\phi_j = \phi_j(\mu)$ depends on $\mu$.

Let us first consider the case $\mu > 0$.
Stripping away positive terms and using $\mu \le 1$ gives
\begin{equation}
 \phi_j \geq \ansatz_j -\cfrac{\dt}{\dx} \psi_{j + 1 / 2}^-
  - \dt \overline{\sig{t} \psi}_j.
\label{eq:phi-lower-bd}
\end{equation}
Next, we want to use our polynomial reconstruction
\eqref{eq:reconstruction-in-cell} and the $\nqs$-point Gauss-Lobatto
quadrature to relate the exact cell mean $\ansatz_j$ to the reconstruction's values at the
cell edge $ \psi_{j + 1 / 2}^-$ and the cell average $\overline{\sig{t} \psi}_j$.
But since only the approximations from the optimization are used to define
the reconstruction, we must first multiply and divide $\ansatz_j$ by
$\psibar_j$. Note that this is possible since $\psibar_j>0$.
Then we use the assumed bound \eqref{eq:ratio-bound}, note carefully that
for $\mu > 0$ indeed $\psi_{j + 1 / 2}^- = \psi_j(x_{j\nqs})$), and apply
$L^\infty$ bounds on the total interaction coefficient.
This gives

\begin{align*}
 \phi_j & \stackrel{\eqref{eq:phi-lower-bd}}{\geq} \frac{\ansatz_j}{\psibar_j} \sum_{i = 1}^\nqs \hat{w}_i \psi_j(x_{ji})
   - \cfrac{\dt}{\dx} \psi_j(x_{j\nqs})
   - \dt \sum_{i = 1}^\nqs \hat{w}_i \sig{t}(x_{ji}) \psi_j(x_{ji}) \\
  & \stackrel{\eqref{eq:ratio-bound}}{\geq} (1 - \veps) \sum_{i = 1}^\nqs \hat{w}_i \psi_j(x_{ji})
   - \cfrac{\dt}{\dx} \psi_j(x_{j\nqs})
   - \dt \sum_{i = 1}^\nqs \hat{w}_i \sig{t}(x_{ji}) \psi_j(x_{ji}) \\
  & \stackrel{\eqref{eq:sigmatmax}}{\geq} \sum_{i = 1}^{\nqs - 1} \hat{w}_i
   \left( 1 - \veps - \dt \sig{t, \max} \right) \psi_{ji}
   + \left( (1 - \veps)\hat{w}_{\nqs} - \cfrac{\dt}{\dx} - \dt \hat{w}_{\nqs}
   \sig{t, \max} \right) \psi_{j \nqs},
\end{align*}

where in the last line we have introduced the notation $\ansatz_{ji} :=
\ansatz_j(x_{ji})$.
One can see that \eqref{eq:CFL} ensures nonnegativity of both terms in the
final expression.
Recalling that $\hat{w}_1 = \hat{w}_\nqs$, the cases $\mu < 0$ and $\mu = 0$ follow
analogously, and
together we have that $\phi_j \ge 0$ at each $\mu$, which shows that
$\ubar^{\timelvl + 1}_j$ is realizable.
\end{proof}

\begin{remark}
In practice, \eqref{eq:ratio-bound} is only approximately enforced by
\eqref{eq:stop}.
However, we have never had problems losing realizability in our numerical
simulations, including those that approach the boundary of realizability.
\end{remark}

\begin{remark}
We have assumed that the reconstruction \eqref{eq:reconstruction-in-cell}
is always used to compute $\overline{\sig{t} \psi}_j$.
However, when $\sig{a}$ and $\sig{s}$ are constant in each cell, one simply
has $\overline{\sig{t} \psi}_j = \sig{t} \ansatz_j$, and therefore the
reconstruction is not needed.
The proof can then be completed with an upper bound on
$\ansatz_j / \psibar_j$, which is similarly easy to enforce approximately,
and one ends up with a slightly less restrictive CFL condition.
However, we did not use this in our implementation.
\end{remark}

The time-step restriction \eqref{eq:CFL} can then be scaled by $\rho$ for the corresponding
time-integration scheme to give realizability of every stage and step in the
scheme.

In practice, to achieve an order $\order$ method for sources $S$
or interaction coefficients $\sig{a}$ or $\sig{s}$ which are not piecewise
degree $\order - 1$ polynomials, one would approximate them using the same
spatial reconstruction techniques that we use for the density to achieve
an order $\order$ approximation of the corresponding terms.
Thus in \thmref{thm:main}, one would not use a value of $\orderS$ larger than
$\order$.

\subsection{Limiting}

The first role of the limiter, then, is first to ensure that the point-wise
values of the polynomial reconstructions $\ansatz_j$ are nonnegative at the
spatial and angular quadrature points.
However, as we will see in \secref{sec:SourceBeam} below, the numerical
solutions using a limiter which only ensures nonnegativity can still contain
spurious oscillations.
Therefore we extend the same limiter to enforce local maximum
principles as well, thereby much more effectively dampening such
oscillations.

\subsubsection{Positivity-preserving limiter}
To preserve nonnegativity we can simply apply a linear scaling limiter.
The limited spatial reconstruction is defined as $\psi_j^\theta(x) :=
(1 - \theta) \psibar_j + \theta \psi_j(x)$; notice that $\theta = 1$
corresponds to no limiting.
For each quadrature point (both in space and angle, though here we suppress
the angular argument) we compute
\begin{align}
\label{eq:Limiter}
 \theta_{ji} = \begin{cases}
   \cfrac{\psibar_j}{\psibar_j-\psi_{ji}}
    & \text{ if } \psi_{ji} < 0, \\
   1 & \text{ else.}
  \end{cases}
\end{align}
Then in each cell we set
$$
\theta = \theta_j := \min_{i=1,\ldots,\nqs} \{ \theta_{ji} \}
$$
(where one should keep in mind that $\theta_j$ still depends on $\mu$).
One immediately sees that this limiter ensures the positivity,
preserves the cell means $\psibar_j$, and following
arguments from \cite{Zhang2010a,Zhang2011}, does not destroy accuracy of
the scheme if $\psibar_j > 0$.

However, it has been remarked in \cite{Zhang2012a} that in some pathological
situations this limiter may reduce the accuracy to second order. See also \cite{AlldredgeSchneider2014} where a similar observation has been made for a realizability limiter in a discontinuous-Galerkin scheme.

\subsubsection{Maximum principle-satisfying limiter}
\label{sec:MaxLimiter}
The limiter we introduce here is a slightly modified version of the
maximum-principle limiter from \cite{Zhang2011a}.

Since we know \emph{a priori} that $\psi$ satisfies a strict maximum
principle $m\leq\psi(x,\mu)\leq M$ for all $x$ and $\mu$, a natural strategy
to dampen artificial oscillations in
numerical solutions is to enforce a local maximum principle.
Specifically, we would like the polynomial reconstruction $\psi_j(x)$ to be
bounded by the data of those cells which influence it.
The corresponding index set of influential nodes is
$\cN_{j, \order} = \{j - \order, \ldots , j + \order \}$ (cf.
\eqref{eq:interpolation}), so the local maximum principle we would like to
enforce is
$$
\min\limits_{\mu, \ell \in \cN_{j,\order}} \psibar_\ell(\mu)
  \leq \psi_j(x, \mu) \leq \max\limits_{\mu, \ell \in
  \cN_{j, \order}} \psibar_\ell(\mu).
$$
However this tends to flatten smooth extrema, so, inspired by the modified
minmod function in \cite{CockburnShuPk}, we relax the strict maximum
principle by setting the maximum principle bounds $M_j$ and
$m_j$ locally as
$$
 M_j := \left(1 + c \frac{\dx}{2} \right)
  \max_{\mu, \ell \in \cN_{j, \order}} \psibar_\ell(\mu)
 \quand
 m_j := \left(1 - c \frac{\dx}{2} \right)
  \min\limits_{\mu, \ell \in \cN_{j, \order}} \psibar_\ell(\mu),
$$
where $c$ is a local bound on the relative derivative of $\psi$, i.e.
$\max_{\mu,x} |\partial_x \psi(x, \mu) / \psi(x, \mu)|$, and
the maximum and minimum in $\mu$ are taken over the angular quadrature
nodes.%
\footnote{
Note that we always want $m_j \ge 0$, so we never take $c > 2 / \dx$.
}
Therefore the maximum principle that we will actually enforce is
$$
m_j \leq \psi_j(x_{ji}, \mu) \leq M_j
$$
for all spatial quadrature points $x_{ji} \in I_j$ and all angular
quadrature points.
To enforce this maximum principle, for each spatial quadrature point we set

\begin{equation}
\label{eq:MPLimiter}
 \theta_{ji} := \begin{cases}
   \cfrac{M_j - \psibar_j}{\psi_{ji} - \psibar_j} & \mbox{if }
    \psi_{ji} > M_j, \\
   \cfrac{m_j - \psibar_j}{\psi_{ji} - \psibar_j} & \mbox{if }
    \psi_{ji} < m_j, \\
   1 & \mbox{otherwise,}
  \end{cases}
\end{equation}
and finally for each cell we choose $\theta_j := \min_i \{ \theta_{ji} \}$.

As with the positivity-preserving limiter, it can be shown
that this limiter does not destroy accuracy \cite{Zhang2010a}.

%% file: Sections/Results.tex
%!TEX root = ../MMEOptimizationKRM.tex

%%%%%%%%%%%%%%%%%%%%%%
\section{Numerical results}
%%%%%%%%%%%%%%%%%%%%%%
\label{sec:Results}
In this section we present results to confirm that our scheme converges
with the expected order, to show the effect of various parameters in the
scheme, and to highlight some of the features of high-order solutions.

Except where otherwise noted, we used the following parameter values:

\begin{tabular}{r @{$\:$} c @{$\:$} l l}
  $\tau$ &=& $10^{-9}$ & Optimization gradient tolerance, \\
  $\veps$ &=& $0.01$ & Optimization tolerance on
   $1 - \ansatz_j / \psibar_j $, \\
  $\{ r_1, \ldots , r_M \}$ &=& $\{ 10^{-8}, 10^{-6}, 10^{-4}\}$
   & Outer regularization loop in optimizer, \\
  $k_r$ &=& $50$ & Number of optimization iterations before \\
   &&& advancing outer regularization loop, \\
  $\nQ$ &=& $40$ & Number of angular quadrature nodes, \\
  $c$ &=& $1$ & Bound on $|\partial_x \psi / \psi|$ in
   maximum-principle limiter.
\end{tabular}

For the angular quadrature we used $(\nqmu / 2)$-point Gauss-Lobatto rules
over both $\mu \in [-1, 0]$ and $\mu \in [0, 1]$.

For the value $q$ determining the number of initialization steps for the
multi-step time integrators we used two for fifth-order simulations and
three for sixth- and seventh-order simulations.

The time step is chosen to fulfill \eqref{eq:CFL} (replacing $\dt$ by
$\dt/\rho$ for the appropriate time integrator) with equality.

In both benchmark problems we use isotropic scattering, $\collision(\psi)
= \frac12 \vint{\psi} - \psi$. 

For the first stage of the first time step, the initial multipliers
for the optimizer are those associated with the normalized isotropic
distribution.
For the following stages and time steps, the initial multipliers are set to
those from the previous stage or step at the same spatial cell.
However, if in these later stages the optimizer cannot converge before
initializing the regularization loop, we first switch back to the multipliers
associated with the normalized isotropic distribution and restart the
optimizer (this time allowing regularization if necessary).
In our experience, the isotropic multipliers are the safest choice for the initial
condition, and therefore this technique reduces the number of times
the regularization must be used.

As in \cite{AllHau12}, if regularization has to be applied to $\ubar_j$, we replace it (in
\eqref{eq:moment-euler-step}) by the regularized moments $\bv(\ubar_j, r)$
\eqref{eq:RegularizedMoments} for which the optimizer converged.
Then \thmref{thm:main} can be applied to ensure that the next iterate will be realizable as well.

In all figures below we plot the zeroth-order component of the
reconstruction $\bu_j(x) = \vint{\basis \psi_j(x, \mu)}$ for $x \in I_j$,
see \eqref{eq:reconstruction-in-cell}.

\subsection{M${}_\nmom$ manufactured solution}

In general analytical solutions for minimum-entropy models are not known.
Therefore, to test the convergence and efficiency of our scheme, we use the
method of manufactured solutions, and we follow the target solution given
in \cite{AlldredgeSchneider2014} but add a spatially and
temporally dependent absorption interaction coefficient.
The solution is defined on the spatial domain $X = (-\pi, \pi)$ with
periodic boundary conditions.

A kinetic density in the form of the entropy ansatz is given by
\begin{subequations}
\label{eq:MFSM3}
\begin{align}
 \phi(t, x, \mu) =& \exp(\alpha_0(t, x) + \alpha_1(t, x) \mu ), \\
 \alpha_0(t, x) =& -K - \sin(x - t) - a,\\
 \alpha_1(t, x) =&  K + \sin(x - t).
\end{align}
\end{subequations}
A source term is defined by applying the transport operator to $\phi$:
$$
\source(t, x, \mu) := \partial_t \phi(t, x, \mu) + \mu \partial_x
\phi(t, x, \mu) + \sig{a}(t, x) \phi(t, x, \mu),
$$
where
$$
 \sig{a}(t, x) := 4 - 4 \cos(x - t)).
$$
Thus by inserting this $\source$ into \eqref{eq:FokkerPlanck1D} (and taking
$\sig{s} = 0$)
we have that $\phi$ is, by construction, a solution of
\eqref{eq:FokkerPlanck1D}.

A straightforward computation shows that $\source \ge 0$ (for
any $a$ or $K$), which means that \thmref{thm:main} will apply to the
resulting moment system.
Furthermore we take
\begin{align*}
 a &= - K + 1 - \log\left( \cfrac{K - 1}{2\sinh(K - 1)} \right)
\end{align*}
so that the maximum value of $\vint{\phi}$ for $(t, x) \in [0, \tf] \times X$
is one.
As $K$ is increased, $\phi$ converges to a Dirac delta at $\mu = 1$.

The moment vector $\bw = \vint{\basis \phi}$ is then a solution of
\eqref{eq:moment-closure}
for $\MN$ models with $\nmom \ge 1$.

We used the final time $\tf = \pi / 5$ and
chose $K = 4$, for which $w_1 / w_0 \in [0.67,0.8]$
(recall that $|w_1 / w_0| < 1$ is necessary for realizability). 
We are using a fairly low value of $K$ because, in order to show convergence
for our scheme with the highest-order ($\order = 7$), we need to be able to
have a tighter control on the errors from the numerical optimization.
When $K$ is higher, these errors are larger and drown out the convergence
in space and time.
In the following, we used the M$_3$ model so that our results included the
effects of the numerical optimization.

We compute errors in the zero-th moment of the solution, which we denote
$w_0(t, x) = \vint{\phi(t, x, \cdot)}$.
Then $L^1$ and $L^\infty$ errors for $u_{0,h}(t, x)$
(that is, the zero-th component of a numerical solution $\U_h$) are
defined as
\begin{align*}
 E^1_h = \int_X \left| w_0(\tf, x) - u_{0,h}(\tf, x) \right| dx
  \quand
 E^\infty_h = \max_{x \in X} \left| w_0(\tf, x) - u_{0,h}(\tf, x) \right|
\end{align*}
respectively. We approximate $u_{0,h}(\tf, x)$ using the same reconstruction
technique as for the scheme for the underlying kinetic density
\eqref{eq:reconstruction-in-cell} and integrate with respect to $\mu$.
Then we approximate the integral in $E^1_h$ using a 100-point Gauss-Lobatto
quadrature rule over each spatial cell $I_j$, and $E^\infty_h$ is
approximated by taking the maximum over these quadrature nodes.

The observed convergence order $\nu$ is defined by
$$
\frac{E^p_{h1}}{E^p_{h2}} = \left( \frac{\dx_2}{\dx_1} \right)^\nu
$$
where for $i \in \{1, 2\}$, $E^p_{hi}$ is the error $E^p_h$ for the
numerical solution using cell size $\dx_i$, for $p \in \{1, \infty\}$.

Convergence tables are given in \tabref{tab:MFSM3} for solutions with a
tighter optimization
tolerance of $\tau = 10^{-11}$.
We observe that the scheme converges with at least its designed order until
the errors are roughly $\cO(\tau)$, where errors from the numerical
optimization halt the convergence.
For many of the solutions on the coarsest grids, the convergence is faster
than designed, likely because the WENO reconstruction is order $2\order-1$ at
the cell interfaces.
The effects are indeed more pronounced for higher orders.

\begin{figure}
\centering
\externaltikz{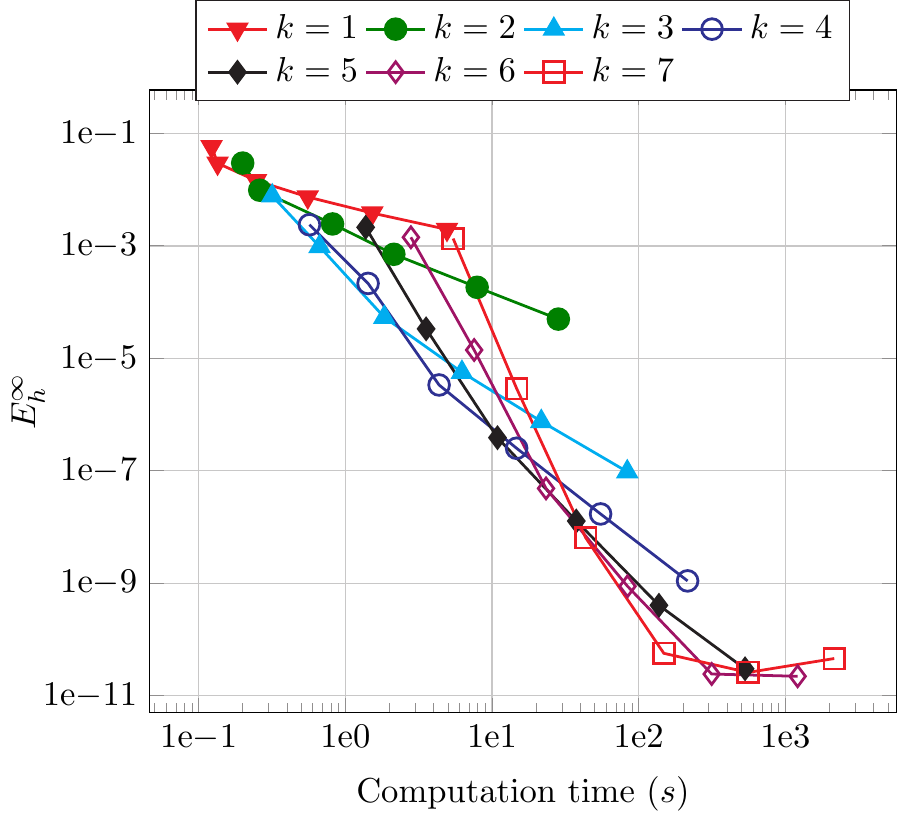}{\input{Images/ConvergenceM3.tex}}
          
\caption{Efficiency of the scheme up to seventh order for the M${}_3$
manufactured solution with optimization tolerance $\tau = 10^{-11}$.}
\label{fig:M3ConvergencePBC}
\end{figure}

In \figref{fig:M3ConvergencePBC} we plot the error of solutions of various
orders against their computation time.
Here we confirm the expectation that for smaller errors, higher-order
solutions require less computation time.

\begin{table}
\centering
\begin{tabular}{r r@{.}l c r@{.}l c r@{.}l c r@{.}l c }
& \multicolumn{3}{c}{$\order = 2 $}& \multicolumn{3}{c}{$\order = 3 $}& \multicolumn{3}{c}{$\order = 5 $}& \multicolumn{3}{c}{$\order = 7 $}\\
\cmidrule(r){2-4} \cmidrule(r){5-7} \cmidrule(r){8-10} \cmidrule(r){11-13} 
$\ngrid $ & \multicolumn{2}{c}{$E^1_h$} & $\nu$ & \multicolumn{2}{c}{$E^1_h$} & $\nu$ & \multicolumn{2}{c}{$E^1_h$} & $\nu$ & \multicolumn{2}{c}{$E^1_h$} & $\nu$\\ \midrule 
 10 & 6 & 532e-02 & ---& 1 & 668e-02 & ---& 4 & 816e-03 & ---& 2 & 696e-03 & ---\\
20 & 1 & 981e-02 & 1.7& 9 & 931e-04 & 4.1& 3 & 622e-05 & 7.1& 5 & 130e-06 & 9.0\\
40 & 3 & 823e-03 & 2.4& 5 & 531e-05 & 4.2& 2 & 517e-07 & 7.2& 3 & 452e-09 & 10.5\\
80 & 1 & 005e-03 & 1.9& 6 & 808e-06 & 3.0& 7 & 546e-09 & 5.1& 5 & 090e-11 & 6.1\\
160 & 2 & 193e-04 & 2.2& 9 & 778e-07 & 2.8& 2 & 427e-10 & 5.0& 4 & 049e-11 & 0.3\\
320 & 5 & 784e-05 & 1.9& 1 & 317e-07 & 2.9& 4 & 636e-11 & 2.4& 5 & 645e-11 & -0.5\\

 & \multicolumn{2}{c}{$E^\infty_h$} & $\nu$ & \multicolumn{2}{c}{$E^\infty_h$} & $\nu$ & \multicolumn{2}{c}{$E^\infty_h$} & $\nu$ & \multicolumn{2}{c}{$E^\infty_h$} & $\nu$\\ \midrule 
 10 & 2 & 963e-02 & ---& 7 & 731e-03 & ---& 2 & 133e-03 & ---& 1 & 347e-03 & ---\\
20 & 9 & 754e-03 & 1.6& 9 & 713e-04 & 3.0& 3 & 343e-05 & 6.0& 2 & 905e-06 & 8.9\\
40 & 2 & 452e-03 & 2.0& 5 & 360e-05 & 4.2& 3 & 879e-07 & 6.4& 6 & 511e-09 & 8.8\\
80 & 7 & 076e-04 & 1.8& 5 & 655e-06 & 3.2& 1 & 280e-08 & 4.9& 5 & 687e-11 & 6.8\\
160 & 1 & 832e-04 & 1.9& 7 & 544e-07 & 2.9& 4 & 063e-10 & 5.0& 2 & 598e-11 & 1.1\\
320 & 4 & 980e-05 & 1.9& 9 & 613e-08 & 3.0& 3 & 017e-11 & 3.8& 4 & 594e-11 & -0.8\\
\end{tabular}
\caption{$L^1$- and $L^\infty$-errors and observed convergence order $\nu$
for the M$_3$ manufactured solution \eqref{eq:MFSM3} with optimization
gradient tolerance $\tau = 10^{-11}$.}
\label{tab:MFSM3}
\end{table}

\subsection{Plane source}
In this test case we start with an isotropic distribution where the initial 
mass is concentrated in the middle of an infinite domain $x \in
(-\infty, \infty)$:
\begin{align*}
 \psiInit(x, \mu) = \psi_{\rm floor} + \delta(x)
\end{align*}
where the small parameter $\psi_{\rm floor} = 0.5 \times 10^{-8}$ is used to
approximate a vacuum.
In practice, a bounded domain must be used, so we choose a domain large
enough that the boundary should have only negligible effects on the
solution:
thus for our final time $\tf = 1$, we take $X = [\xL, \xR] = [-1.2, 1.2]$.
At the boundary we set
\begin{align*}
 \psiL(t, \mu) \equiv \psi_{\rm floor} \quand
 \psiR(t, \mu) \equiv \psi_{\rm floor}
\end{align*}
We set $\sig{s} = 1$ and $\sig{a} = 0$.

All solutions here are computed with an even number of cells, so the delta
function in the initial condition lies on a cell boundary.
Therefore we approximate the delta function by splitting
it into the cells immediately to the left and right.

Since this is a highly non-smooth problem, the choice of $c$ is not
immediately obvious.
Indeed, the solution contains such strong gradients, so choosing
$c$ too small can greatly reduce the accuracy of the simulation.
After some numerical experimentation, we observed that in the $1000$-cell
seventh-order simulations the value $c = 15$
generally gave a good trade-off between smoothness in the solution without
adding too much diffusivity.
In the $4000$-cell fourth-order simulations, where the non-smoothness is
further resolved, $c = 10$ gave the best results.

\begin{figure}[h!]
\centering
\subfloat[M$_3$ with the positivity-preserving limiter]
{\externaltikz{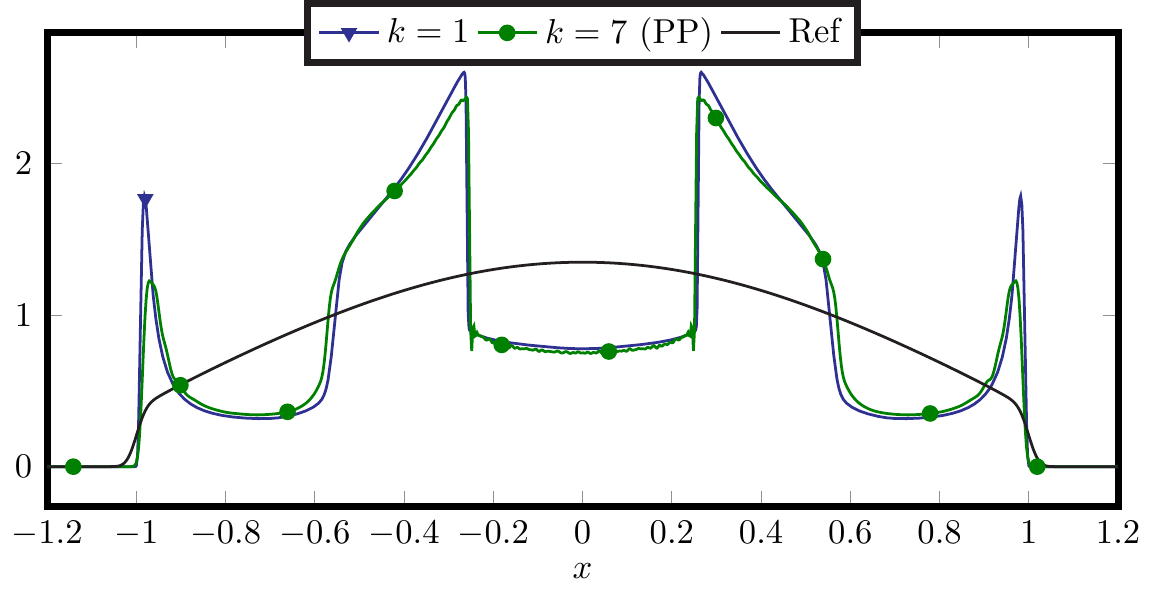}
   {\input{Images/PlanesourceM3a.tex}
    \label{fig:plane-M3-PP}}}\\
\subfloat[M$_3$, $c = 15$]{\externaltikz{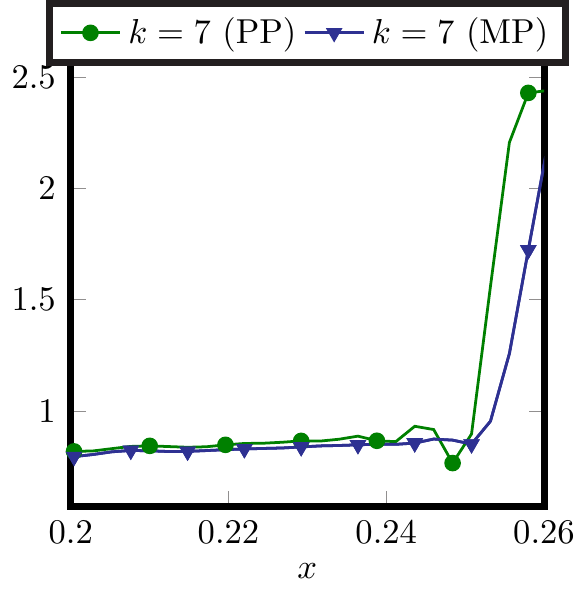}
   {\input{Images/PlanesourceM3b.tex}
    \label{fig:plane-M3-PP-vs-MP-mid}}}
\subfloat[M$_3$, $c = 15$]{\externaltikz{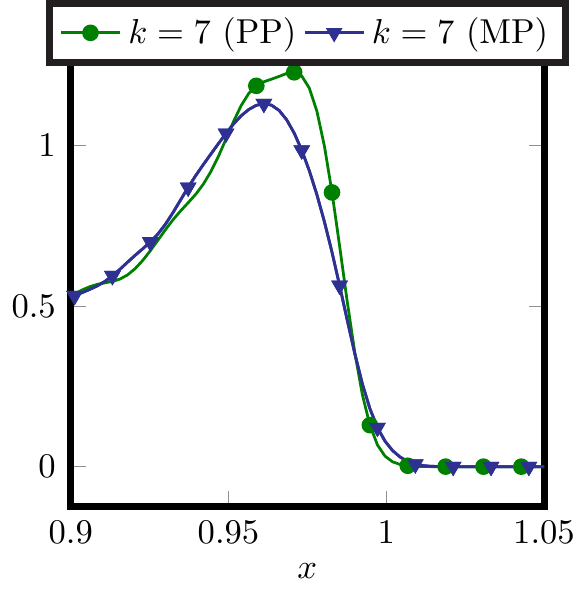}
   {\input{Images/PlanesourceM3c.tex}
    \label{fig:plane-M3-PP-vs-MP-front}}}\\
 \caption{The local particle density $u_0$ from first- and
 seventh-order solutions to the plane-source problem with the M$_3$ $t = 1$.
 First-order solutions were computed with $\ngrid = 10000$ cells and
 seventh-order with $\ngrid = 1000$ cells.}
 \label{fig:plane-M3}
\end{figure}

\begin{figure}
\centering  
\subfloat[M$_7$, $\order = 4$, $c = 10$]{\externaltikz{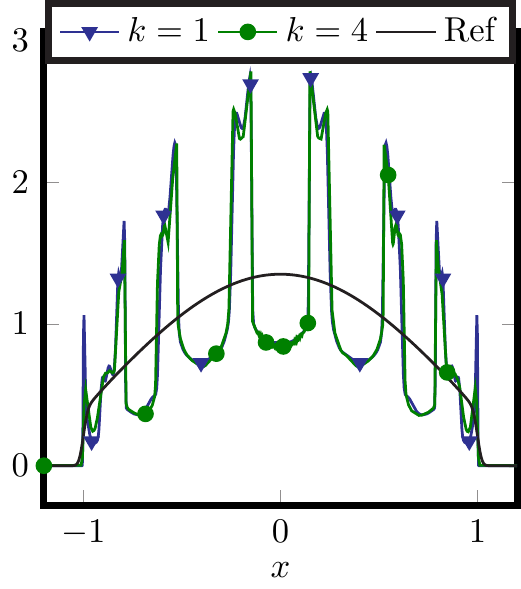}
   {\input{Images/PlanesourceM7a.tex}
    \label{fig:plane-M7-k4}}}
\subfloat[M$_7$, $\order = 7$, $c = 15$]{\externaltikz{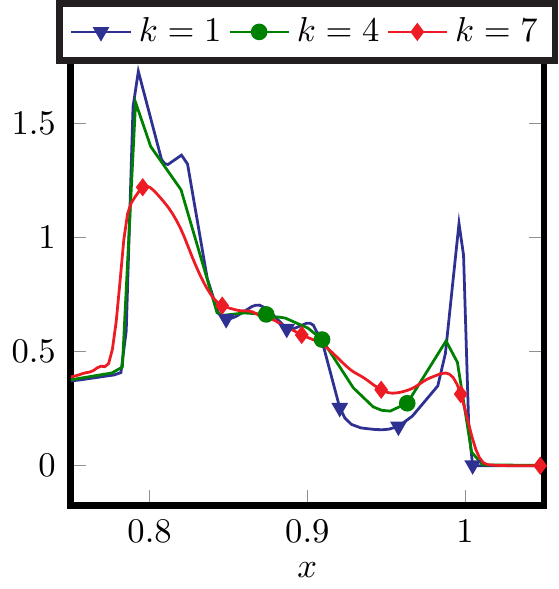}
   {\input{Images/PlanesourceM7b.tex}
    \label{fig:plane-M7-k7}}}\\
\caption{The local particle density $u_0$ from first-, fourth- and
seventh-order solutions to the plane-source problem with the
M$_7$ model at $t = 1$.
First-order solutions were computed with $\ngrid = 10^4$ cells;
fourth-order with $\ngrid = 4000$ cells; and
seventh-order with $\ngrid = 1000$ cells.}
\label{fig:plane-M7}
\end{figure}

In \figref{fig:plane-M3} and \figref{fig:plane-M7}, we plot several solutions
at the final time for the M$_3$ and M$_7$ models, respectively, as well as a
reference solution%
\footnote{
The reference solution was computed using a first-order P$_{199}$ method on
a grid with $\ngrid=4000$ cells.
}
of the kinetic equation \eqref{eq:FokkerPlanck1D}.
We consider numerical solutions of orders $\order \in \{1, 4, 7\}$.
The first-order solutions are included to indicate our best guess of the true
$\MN$ solutions, and they largely agree with those presented in
\cite{Hauck2010} (although we have computed solutions on a much finer grid
in order to better resolve sharp peaks in the solutions).

In \figref{fig:plane-M3}, we present seventh-order M$_3$ solutions.
In \figref{fig:plane-M3-PP}, we see that while the seventh-order solution
using only the positivity-preserving limiter closely matches the
highly-resolved first-order solution, some spurious oscillations around
$x = \pm 0.25$ remain.
These oscillations are largely removed by using the maximum-principle
limiter with $c = 15$, as we see in \figref{fig:plane-M3-PP-vs-MP-mid},
though in \figref{fig:plane-M3-PP-vs-MP-front}
we see that this solution is also somewhat more diffusive.

The M$_7$ solution to the plane-source problem is more oscillatory.
\figref{fig:plane-M7-k4} shows that our kinetic scheme with fourth-order
reconstructions and the maximum-principle limiter performs well.
However, seventh-order solutions are notably more diffusive, as shown at the
leading edge of the solution in \figref{fig:plane-M7-k7}.

That the high-order solutions are more diffusive than the true $\MN$
solutions is not necessarily a disadvantage: the more diffusive solutions
are actually closer to the reference kinetic solution.
Indeed, the $\MN$ equations are a spectral method for the original kinetic
equation, and thus should not be applied to a non-smooth problem like this
one without filtering to avoid the Gibbs phenomenon.
Exactly how to apply such a filter for moment models kinetic equations is
a topic of ongoing research \cite{McClarren20105597,Radice2013648}, but here
it seems that the
maximum-principle limiter already filters the solution somewhat.

\subsection{Source-beam}
\label{sec:SourceBeam}
Finally we present a discontinuous version of the source-beam problem from
\cite{Frank:1471763}.
The spatial domain is $X = [0,3]$, and
\begin{gather*}
 \sig{a}(x) = \begin{cases}
   1 & \text{ if } x\leq 2\\
   0 & \text{ else}
  \end{cases}, \quad
 \sig{s}(x) = \begin{cases}
   0 & \text{ if } x\leq 1\\
   2 & \text{ if } 1<x\leq 2\\
   10 & \text{ else}
  \end{cases}, \quad
 \source(x) = \begin{cases}
   1 & \text{ if } 1\leq x\leq 1.5\\
   0 & \text{ else}
  \end{cases},
\end{gather*}
with initial and boundary conditions
\begin{gather*}
 \psiInit(x, \mu) \equiv \psi_{\rm floor}, \\
 \psiL(t, \mu) = \beta \exp(-\gamma(\mu-1)^2)
 \quand
 \psiR(t, \mu) \equiv \psi_{\rm floor},
\end{gather*}
for $\gamma = 10^5$, normalization constant
$\beta = \vint{\exp(-\gamma(\mu - 1)^2)}^{-1}$, and
$\psi_{\rm floor} = 0.5 \times 10^{-10}$.
 
\begin{figure}[h!]
\setlength\figureheight{0.25\textheight}
\setlength\figurewidth{0.5\textwidth}
\centering
\subfloat[Seventh-order $\MN$ solutions, $\ngrid = 150$ \label{fig:SourceBeamMPAll}]{\externaltikz{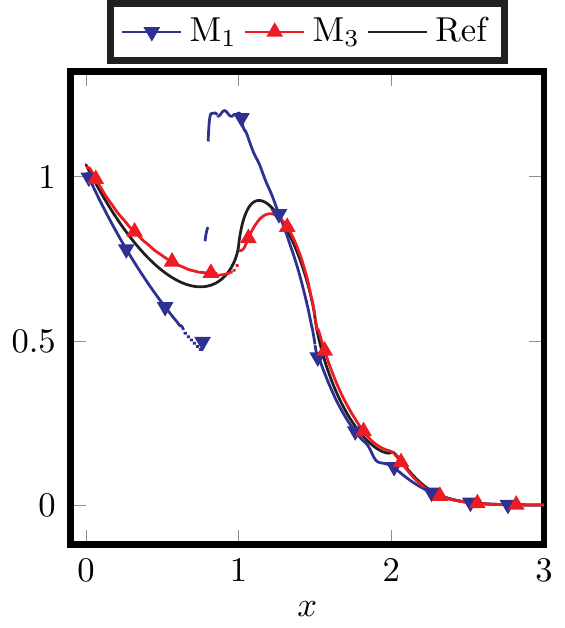}{\input{Images/SourceBeamRef.tex}}}
\subfloat[Order comparison for M$_1$ \label{fig:SourceBeamOrder}]{\externaltikz{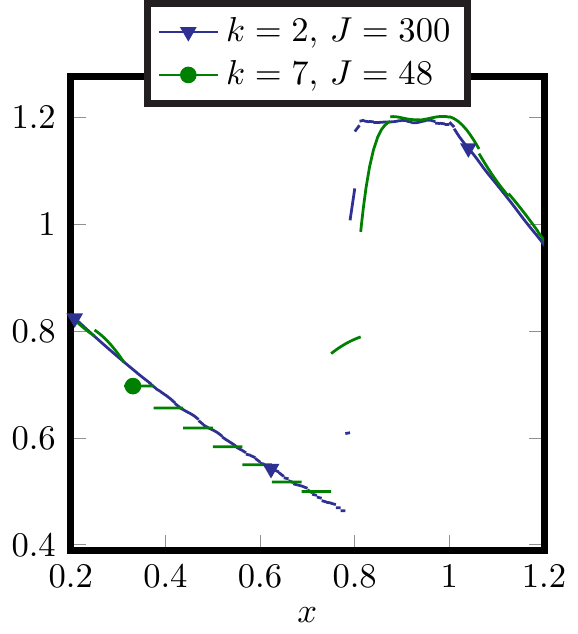}{\input{Images/SourceBeamEfficiency.tex}}}
          
\caption{Seventh-order approximation of several $\MN$ solutions for the source-beam testcase and comparison of the seventh-order approximation on a coarse grid with the second-order approximation on a finer grid.}

\end{figure}

$\MN$ solutions for this problem are shown in \figref{fig:SourceBeamMPAll}
using $\ngrid = 150$ cells and seventh-order reconstructions, using the maximum-principle limiter with $c=1$, along with a
reference solution.%
\footnote{
The reference solution was computed using a first-order P$_{99}$ method on
a grid with $\ngrid=2000$ cells.
}
We see that increasing the moment order to $\nmom = 3$ qualitatively
improves the solution significantly.

\figref{fig:SourceBeamOrder} shows some strengths and weaknesses of our
scheme for the M$_1$ model, where for this problem the shock is the
strongest.
When compared to a much more finely resolved low-order approximation, the
seventh-order approximation on a coarse grid nicely fits the features of the
solution despite the discontinuous physical parameters.
However, the incoming beam can only be resolved up to first order (note the
piecewise-constant approximations behind the beam), resulting in a more
diffusive solution at the shock around $x=1$.

\begin{figure}
      \setlength\figureheight{0.25\textheight}
      \setlength\figurewidth{0.5\textwidth}
  \centering  
 \subfloat[A closer look at oscillations in the M$_1$ solution for different
 limiter configurations.]{
\externaltikz{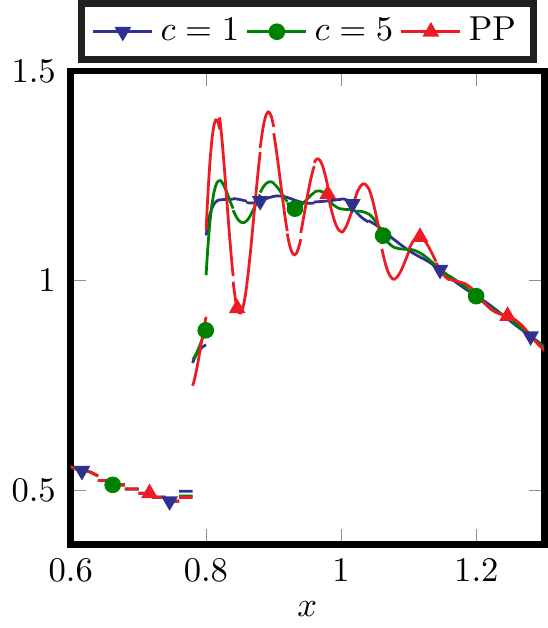}{\input{Images/SourceBeam2.tex}}
\label{fig:source-beam-M1}
  		} 
 \subfloat[Examining the diffusivity in the M$_2$ model for different limiter
  configurations.]
  {\externaltikz{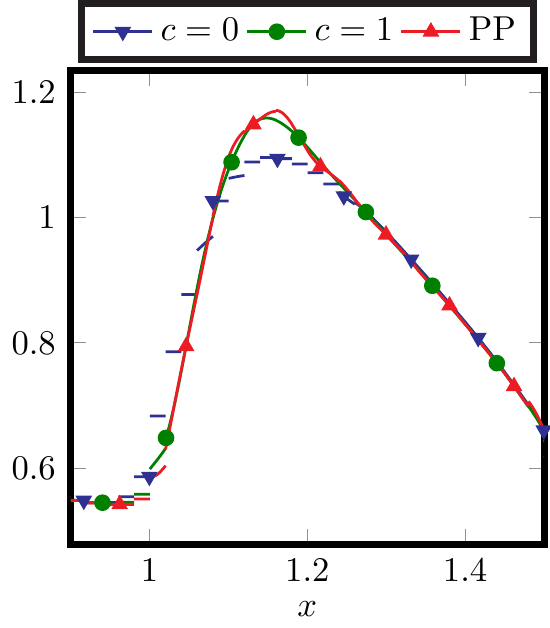}
  {\input{Images/SourceBeam3.tex}}
  \label{fig:source-beam-M2}}

 \caption{The effects of different limiting in the source-beam test at
  $t = 2.5$ for seventh-order solutions with $\ngrid = 150$ cells.}
\label{fig:source-beam-limiting}
\end{figure}

The source-beam problem is particularly well-suited to test the
oscillation-dampening effects of the limiters.
In \figref{fig:source-beam-limiting} we compare several seventh-order
solutions.
First in \figref{fig:source-beam-M1} we see in the M$_1$ solution that the
positivity-preserving limiter does not
dampen spurious oscillations while the maximum principle-preserving limiter
does a better job even for a relatively large value of $c$.
Second, in \figref{fig:source-beam-M2}, we show with the M$_2$ model that for
small values of $c$, the maximum principle-preserving limiter is too
diffusive.

%% file: Sections/Conclusions.tex
%!TEX root = ../MMEOptimizationKRM.tex

%%%%%%%%%%%%%%%%%%%%%%
\section{Conclusions and outlook}
%%%%%%%%%%%%%%%%%%%%%%
\label{sec:Conclusions}

In this paper we describe how to implement a kinetic scheme of (in principle)
arbitrarily high order for entropy-based moment closures of linear kinetic
equations in one space dimension.
For spatial reconstructions we use the well-known WENO method to reconstruct
the underlying entropy ans\"atze using interpolating polynomials,
and time integration is performed using multi-step SSP methods.
These SSP time integrators play a key role in allowing us to give a time-step
restriction which guarantees that the moments stay in the realizable set.
The other key component is a limiter, which not only ensures
positivity of the polynomial reconstructions on a spatial quadrature set,
but also enforces a local maximum principle which dampens spurious oscillations
in numerical solutions.

We performed convergence tests with a manufactured solution that included
the effects of a space- and time-dependent absorption interaction coefficient,
and these results validated that the scheme is converging at least as fast as
expected, and often faster at lower resolutions with higher orders.
The convergence tests also showed that errors from the
numerical optimization routine needed for the angular reconstructions limit
the overall accuracy of the scheme.
This indeed eliminates the benefit of going beyond a certain order (depending
on the optimization tolerance $\tau$).
However,
our manufactured-solution tests showed that increases in efficiency can be
obtained before the optimization errors dominate the solution.

Using the plane-source, an challenging highly non-smooth benchmark problem,
we showed that with the maximum-principle limiter accurate solutions can be
obtained which even limit some of the spurious oscillations due to Gibbs
phenomena in the true $\MN$ solutions, thus pushing our solutions closer to
the kinetic solution.
However, the choice of the parameter $c$ plays an important role in the
approximation quality of the scheme, and a good value is not always
available \emph{a priori}.
With another benchmark problem, the source beam, we also demonstrated the 
benefit of using the maximum principle-preserving limiter.
Again, one must strike a balance between flattening smooth extrema and
dampening spurious oscillations with the choice of its parameter.

Compared to the discontinuous-Galerkin implementation in
\cite{AlldredgeSchneider2014}, where realizability preservation is complex due
to the structure of the set of realizable moments,
realizability preservation on the kinetic level is much simpler.
However, the wider stencils in the WENO reconstruction process will
influence the overall parallelizability of the scheme.

Future work should continue to work toward practical implementations of
entropy-based moment closures.
Models in two and three spatial dimensions should be implemented, and here a
notable challenge is the increasing number of angular quadrature points
that will be needed.
Indeed, our reconstructions are performed at every angular quadrature point,
so more efficient WENO techniques will be necessary.
Other collision models should also be considered.
The Laplace-Beltrami operator in the Fokker-Planck equation does not fall
under the types of collision operators considered here but is an important
model for problems with forward-peaked scattering.
This appears, for example, in important applications such as radiotherapy.
Finally, since we have only considered explicit time-stepping schemes,
our time-step restriction scales with the mean-free path.
Implicit-explicit or asymptotic-preserving schemes should be developed to
handle moment models near diffusive or fluid regimes without requiring
extremely small time steps.